\documentclass[11pt]{article}
\usepackage[utf8]{inputenc}
\usepackage[T1]{fontenc}
\usepackage[toc,page]{appendix}
\usepackage{newtxtext, newtxmath}%--- fonts

\usepackage{amsthm, amsmath, amsfonts} %ams packages
\usepackage{faktor}
\usepackage{thmtools, thm-restate}
\usepackage[hypcap=false]{caption}
\usepackage{tabu}% for tables
\usepackage[shortlabels]{enumitem} %lists
\usepackage[cal=boondoxo]{mathalfa}%mathcal
\usepackage{indentfirst} %indents first paragraph
\usepackage{xspace} %handle space after macros
\usepackage{todonotes}
\makeatother
\usepackage{titlesec} % custom section styles
\titleformat{\section}{\large\bfseries\scshape}{\thesection.}{0.5em}{}[{\vspace{0.5em}\titlerule \vspace{0.5em}}]
\titleformat*{\subsection}{\large\scshape}
\declaretheoremstyle[
    spacebelow=6pt, 
    headfont=\normalfont\bfseries, 
    bodyfont = \normalfont,
    postheadspace=1em, 
    qed=$\blacksquare$, 
    headpunct={:}]{boldstyle}

\declaretheorem[name=Theorem, parent=section]{theorem}
\declaretheorem[name=Conjecture, sibling=theorem]{conjecture}
\declaretheorem[name=Corollary, sibling=theorem]{corollary}
\declaretheorem[name=Lemma, sibling=theorem]{lemma}
\declaretheorem[name=Proposition, sibling=theorem]{proposition}
\declaretheorem[name=Definition, style=definition, sibling=theorem]{definition}
\declaretheorem[name=Problem, style=definition, , sibling=theorem]{problem}
\declaretheorem[name=Remark, style=remark, sibling=theorem]{remark}
\declaretheorem[name=Example, style=definition, sibling=theorem]{example}
\declaretheorem[name=Step, style=definition, numbered=no]{step}

\renewcommand{\cal}[1]{\mathcal{#1}}
\newcommand{\Mod}[1]{\ (\mathrm{mod}\ #1)}

\def\N{\ensuremath{\mathbb N}\xspace} %natural numbers
\def\R{\ensuremath{\mathbb R}\xspace} %real numbers
\def\Q{\ensuremath{\mathbb Q}\xspace} %rationals
\def\Z{\ensuremath{\mathbb Z}\xspace} %integers
\def\P{\ensuremath{\mathfrak p}\xspace} %pattern
\def\F{\ensuremath{\mathbb F}\xspace} %field
\def\OK{\ensuremath{\mathcal{O}_K}\xspace} %O_K
\def\dim{\ensuremath{\mathrm{dim}}\xspace} %diameter
\def\f{\ensuremath{\mathcal f}\xspace}
\def\g{\ensuremath{\mathcal g}\xspace}
\def\h{\ensuremath{\mathcal h}\xspace}
\usepackage{hyperref} %hyperref

\title{Pattern Problems related to the Arithmetic Kakeya Conjecture}
\author{Charlie Cowen-Breen, Elene Karangozishvili, Narmada Varadarajan and Thomas Wang}
\date{}

\begin{document}

\maketitle
\begin{abstract}
We study a variety of problems about homothets of sets related to the Kakeya conjecture.
In particular, we show many of these problems are equivalent to the arithmetic Kakeya conjecture of Katz and Tao.
We also provide a proof that the arithmetic Kakeya conjecture implies the Kakeya conjecture for packing dimension, as this implication was previously only known for Minkowski dimension.

We consider several questions analogous to the classical results of Stein and Bourgain about the Lebesgue measure of a set containing a sphere centered at every point of $[0,1]^n$, where we replace spheres by arbitrary polytopes.
We give a lower bound for polytopes in $\R^n$, and show that this is sharp for simplices.
Finally, we generalize number theoretic methods of Green and Ruzsa to study patterns in number fields and thereby provide upper bounds for several of these homothet problems.
\end{abstract}

\section*{Acknowledgements}
We thank Tam\'as Keleti for his invaluable guidance advising this project as part of the Budapest Semesters in Mathematics summer research experience.
\section{Introduction}
\subsection{Geometric motivation}

Our journey begins with a classical and deep result---discovered first by Stein for $n\geq3$, and then independently by Bourgain \cite{Bourgain} and Marstrand \cite{Marstrand} for $n=2$---that any set $B \subset \R^n$ containing a sphere $S^{n-1}$ centered at each point of $\R^n$ must have positive Lebesgue measure.
This was sharpened by Wolff \cite{Wolff1}, \cite{Wolff2} and Mitsis \cite{Mitsis} to show the following. 
Suppose that $B$ contains a sphere centered at each point of a set $S \subset \R^n$. 
If $\dim_H S\leq 1$, then $\dim_H B\geq n-1+\dim_H S$.
If $\dim_H S > 1$, then $B$ must have positive Lebesgue measure, hence $\dim_H B = n$.
Moreover, Talagrand \cite{Talagrand} and Mitsis \cite{Mitsis} showed that this is sharp for the threshold $\dim_H S = 1$: there exists a set of Lebesgue measure zero containing a sphere centered at each point of a straight line $S$ in $\R^n$.

 Keleti, Nagy, and Shmerkin showed, perhaps surprisingly, that when $n=2$, these statements do not hold if the circles are replaced by the boundaries of axis-parallel squares. 
 Namely, for every $s\in[0,2]$, there exist compact sets $S,B\subset\R^2$ such that $B$ contains the boundary of a square centered at each point of $S$, $\dim_H S = s$, and $\dim_H B = 1$ (\cite{Squares}). 
 Moreover, this result holds identically for all polygons and, indeed, generalizes for all polytopes in $\R^n$ due to Chang, Cs{\"o}rnyei, Hera, and Keleti \cite{Chang}.
 
What does not hold identically for all polygons and polytopes is the corresponding question for Minkowski
and packing dimensions. We will illustrate this by contrasting the results for cubes and simplices.
In \cite{Cubes}, Thornton generalized a result of Keleti et al. to show that if $B \subset \R^n$ contains the boundary of an axis-parallel $n$-cube centered at each point in $[0,1]^n$, then for $\dim = \dim_M$ or $\dim_P$, $\dim B\geq n-\frac{1}{2n}$, and this is sharp.

For simplices, on the other hand, we show here that if $B\subset \R$ contains a homothetic copy of an $n$-simplex centered at each point in $\R^n$, then $\dim B \geq n - 1/(n + 1)$, and this is sharp as well. In fact, we show that this lower bound holds for all convex polytopes in $\R^n,$ although it is likely not sharp.

The general geometric problem is of the following type.
\begin{problem}
Let $\P \subset \R^n$ be a set not containing its center.
If $B \subset \R^n$ contains a homothet of $\P$ centered at each point of $[0,1]^n$, what is the dimension of $B$?
\end{problem}

We are most interested in the case when $\P$ is the boundary of a closed convex body.
For convenience, for the rest of this paper we will write ``convex polytope'' instead of ``boundary of a convex polytope'', i.e.\ we say $\P$ is an $n$-simplex if $\P$ is the boundary of an $n$-simplex.
We summarize the results in the following table.
Let $\g_n^M(\P)$ (resp.\ $\g_n^P(P)$ and $\g_{n}^{H}(\P)$) denote the minimal Minkowski (resp.\ packing and Hausdorff) dimension of a set containing a homothet of $\P$ centered at each point of $[0,1]^n$.
When the dimension $n$ is clear from context, we omit the subscript.
%\todo{Citation for $\dim_H(B)$???}
\begin{center}
\bgroup
\def\arraystretch{1.5}
\begin{tabular}{|c|c|c|c|c|}
\hline
$\P$& $\g^H(\P)$ & $\g^P(\P)$ or $\g^{M}(\P)$\\
\hline
$n$-simplex & $n-1$ & $ n - \frac{1}{n+1}$ \small{(Theorem \ref{thm:simplex})}  \\
\hline
$n$-cube & $n-1$ & $n - \frac{1}{2n}$ \small{\cite{Cubes}}  \\
\hline
convex polytope with $m$ faces & $n-1$ & $n - \frac{1}{n+1} \leq \g \leq n - \frac{1}{m}$ \small{(Theorems \ref{thm:polytopes} and \ref{thm:simplex})}  \\
\hline
$S^{n-1}$ & $n$ & $n$ \small{(\cite{Bourgain, Marstrand, Stein})} \\
\hline
\end{tabular}
\egroup
\captionof{table}{The minimal dimension of a set containing a homothet of $\P$ centered at each point of $[0,1]^n$.}\label{table:geom}
\end{center}
Based on these results, one might expect that for a sequence $(\P_k)_{k\in\mathbb{N}}$ of polytopes which converge to the sphere $S^{n-1}$ (in the Hausdorff metric), the minimal Minkowski or packing dimension of a set $B\subset\R^n$ containing a homothetic copy of $\P_k$ centered at each point in $\R^n$ converges to $n$ as $k\to\infty$.

In this paper we show that for some choices of $(\P_k)_{k \in \N}$, this statement is implied by the arithmetic Kakeya conjecture of Katz and Tao, which we discuss next. Moreover, if we consider $S$ to be a straight line instead of all of $\R^n$, as Talagrand and Mitsis did, then the corresponding statement is equivalent to the arithmetic Kakeya conjecture.

\subsection{The Kakeya conjecture}

Each of the problems in the previous section could rightly be classified as a variant of the celebrated Kakeya conjecture. Here we state two versions of the Kakeya conjecture; a more complete introduction to the Kakeya conjecture can be found in \cite{Kakeya}.
\begin{conjecture}[Kakeya conjecture]\label{conj:KC}
A Besicovitch set is a subset of $\R^n$ containing a unit line segment in every direction.
If $B \subset \R^n$ is a Besicovitch set, then
\begin{enumerate}[(a)]
    \item $\dim_H B = n$.
    \item $\dim_P B = n$.
    \item $\dim_M B = n$.
\end{enumerate}
\end{conjecture}
Clearly (a) $\implies$ (b) $\implies$ (c), but none of the converse implications is known.
First formulated by Katz and Tao \cite{KatzTao}, the Arithmetic Kakeya conjecture (AKC) is an additive-combinatorial statement intended to imply Conjecture \ref{conj:KC}(c). 
This proof of this implication dates back to a ``slicing'' argument of Bourgain \cite{BourgainKakeya}.
We will show something stronger, namely
\begin{theorem}\label{thm:AKCimpliespackingKC}
The Arithmetic Kakeya conjecture implies the Kakeya conjecture for packing dimension .
\end{theorem}

 The following are two equivalent combinatorial statements of AKC, as listed in \cite{GreenRuzsa}. 
\begin{conjecture}[Arithmetic Kakeya conjecture]\label{conj:AKC}
(\cite{GreenRuzsa})
Let $k$ and $N$ be positive integers.
\begin{enumerate}[(i)]
    \item Let $F_k(N)$ be the size of the smallest set of integers that contains an $k$-term arithmetic progression with common difference $d$, for each $d \in [N]$.
    Then,
\[
\lim_{k \to \infty} \lim_{N \to \infty} \frac{\log F_{k}(N)}{\log N} = 1.
\]
    \item Let $F_k'(N)$ be the size of the smallest set of integers that contains an $k$-term arithmetic progression with common difference $d$, for $N$ distinct integer values of $d$.
    Then,
\[
\lim_{k \to \infty} \lim_{N \to \infty} \frac{\log F_{k}'(N)}{\log N} = 1.
\]
\end{enumerate}
\end{conjecture}

\subsection{Pattern problems}
The main objective of this paper is to establish a connection between AKC and the sort of problems considered in section 1.1.
Broadly defined, our setup is as follows.
Let 
\[
\cal{P} = \{ \P_{x,r}: x \in X, r \in R\}
\]
be a collection of subsets $\P_{x,r} \subset M$ of a (commutative) ring $M$, where $X$ and $R$ are indexing sets.
$M$ will typically also be a module over a suitable commutative ring $R$. (The canonical example is the structure of $\R^n$ as an $\R$-module.)

\begin{problem}[Pattern problem]
Let $\cal{P}$ be a family as above.
\begin{enumerate}[(a)]
    \item If $B$ contains a set $\P_{x,r}$ for every $x$ from some subset $S \subset X$, what can we say about the size of $B$?
    \item If $B$ contains a set $\P_{x,r}$ for every $r$ from some subset $S \subset R$, what can we say about the size of $B$?
\end{enumerate}
\end{problem}
The notion of size may refer to cardinality, measure, or dimension -- whichever yields the most nontrivial information about $B$.
Of course, this problem is only interesting when the sets $\P_{x,r}$ are themselves of interest.

\begin{definition}
Where addition and multiplication are suitably defined, we say the set $x + r \cdot \P$ is a \textit{$\P$-pattern} with \textit{basepoint} $x$ and \textit{scaling factor} $r$.
\end{definition}
For example, a $k$-term arithmetic progression $a+d, \dots, a+kd$ is a $[k]$-pattern with basepoint $a\in \Z$ and scaling factor $d\in \Z$.\footnote{Note that the basepoint of a progression is not a member of it.}
Similarly, if $\P$ is a convex polytope in $\R^n$, a homothet $x + r \cdot \P$ with $x \in \R^n$ and $r \in \R$ is a $\P$-pattern.
Note that $x$ and $r$ belong to different sets, and here we exploit the $\R$-module structure of $\R^n$ to define ``multiplication''.

The Kakeya conjecture is a problem of type (b): if $B \subset \R^n$ contains a $[0,1]$-pattern (a unit line segment) for every scaling factor in $S^{n-1}$ (direction), what can we say about the size of $B$?
The pattern problem has a rich history --- see the survey paper \cite{survey} by Keleti.

For a fixed family $\cal{P}$ we say the problems of type (a) and (b) are \textit{dual} to each other.
For example, as AKC is a type (b) problem, we consider its dual problem of type (a) as stated below, which we will show in \autoref{sec:finitepatterns} to be equivalent to AKC. 
In other words, in the limit, demanding that ``for every scaling factor, there exists a basepoint'' is equivalent to demanding that ``for every basepoint, there exists a scaling factor''.

\begin{conjecture}[Arithmetic pattern conjecture]\label{conj:AMpatterns}
Let $k$ and $N$ be positive integers.
\begin{enumerate}[(i)]
    \item Let $G_{k}(N)$ be the size of the smallest set of integers that contains a $k$-term arithmetic progression with basepoint $x$, for each $x\in[N]$. Then,
\[
\lim_{k\to\infty} \lim_{N\to\infty} \frac{\log G_{k}(N)}{\log N} = 1.
\]
    \item Let $G'_{k}(N)$ be the size of the smallest set of integers that contains a $k$-term arithmetic progression with basepoint $x$, for $N$ distinct integer values of $x$. Then,
\[
\lim_{k\to\infty} \lim_{N\to\infty} \frac{\log G'_{k}(N)}{\log N} = 1.
\]
\end{enumerate}
\end{conjecture}

Note that not all dual pattern problems are equivalent. 
For instance, although a set containing a circle with every center in $\R^2$ must have positive Lebesgue measure (as discussed earlier), there are planar sets of measure zero containing a circle with every radius $r>0$, as seen in \cite{ThinCircles}, \cite{ThinCircles2}.

%For consistency, we will use the letter $f$ for problems of type (a1) and (b1), and $g$ for those of type (a2) and (b2).
%When the notion of size is cardinality, we use the symbol $F$ (resp.\ $G$), and when it is dimension, we use the symbol $\mathcal{f}$ (resp.\ $\g$).

\subsection{Statement of main results}

In \autoref{sec:finitepatterns}, we study finite patterns over the integers, such as AKC and its dual.
We define a larger class of pattern conjectures and show that all of them imply AKC, and the converse implication holds for several of them. 
For example, define $1/[k] := \{1/i: i \in [k] \}$, and $G'_{1/[k]}(N)$ to be the size of the smallest set of rational numbers that contains a $1/[k]$-pattern with basepoint $a$, for $N$ distinct integer values of $a$.
\begin{theorem}\label{thm:3patterns}
With the same notation as above,
\[
\lim_{k \to \infty} \lim_{N \to \infty} \frac{\log G'_{1/[k]}(N)}{\log N} = \lim_{k \to \infty} \lim_{N \to \infty} \frac{\log G'_{k}(N)}{\log N} = \lim_{k \to \infty} \lim_{N \to \infty} \frac{\log F'_{k}(N)}{\log N}.
\]
\end{theorem}
Additionally, we give a construction to upper bound $G_k'(N)$, complemented by lower bounds for some small values of $k$.

In \autoref{sec:infinitepatterns}, we study the case when Minkowski and packing dimensions are the correct notions of size.
In these problems, either the set of basepoints or the set of scaling factors is usually infinite.
Let $\f^{M}(k)$ and $\f^{P}(k)$ be the infimum of the Minkowski (resp.\ packing) dimension of a set in $\R$ containing a $k$-term arithmetic progression with every common difference in $[0,1]$.
As a corollary of a more general statement in higher dimensions, it will follow that
\begin{theorem}\label{thm:1dAKCequivs}
With the same notation as above,
\[
\lim_{N \to \infty} \frac{\log F'_{k}(N)}{\log N} = \f^{M}(k) = \f^{P}(k).
\]
\end{theorem}
Finally, through some gymnastics involving intermediate equivalences, we prove \autoref{thm:AKCimpliespackingKC}.

In \autoref{sec:geometry}, we consider problems in higher-dimensional Euclidean space.
We state and prove the results listed in \autoref{table:geom}; we provide bounds for the minimal dimension of a set $B \subset \R^n$ containing a homothet of a convex polytope $\P$ centered at each point of $[0,1]^n$.
One might expect a continuity-type result, that if a sequence of polytopes ``converges'' to the sphere, then the minimal dimension of a corresponding set $B$ converges to $n$.
We show that this is essentially as hard as AKC with the following result.
%Let $\mathcal{g}_n^{M}(\P)$ and $\g_n^{P}(\P)$ denote the minimal dimension (for the cases we consider, Minkowski and packing are the same) of a subset of $\R^n$ containing a scaled copy of $P\subset\R^n$ based at every point of the line segment $S=[0,1]\times\mathbf{0}^{n-1}$. Then
\begin{theorem}\label{thm:harmonicpolytope}
For any $n>1$, there is a sequence of convex polytopes $(\P_k)_{k\in\mathbb{N}}$ which converge in the Hausdorff metric to the sphere $S^{n-1}$, such that AKC is equivalent to the statement that $\g_n(\P_k) \to n$ as $k\to\infty$, where $\g_n$ denotes either $\g_n^P$ or $\g_n^M$.
 %    Moreover,\[  \mathcal{g}_n(P_k) = n-1 + \lim_{N\to\infty} \frac{\log F_k(N)}{\log N} \]
\end{theorem}
Just as with our earlier proofs, we will actually prove something stronger by establishing a relationship between $\g_n^{P}(\P_k)$, $\g_n^{M}(\P_k)$ and $F_k(N)$.
%In light of the results of Wolff and Mitsis discussed in section 1.1, AKC can thus be seen as the continuity statement that
%\[ \lim_{k\to\infty} \mathcal{g}_n(P_k) = \mathcal{g}_n(\lim_{k\to\infty} P_k). \]
%In particular, if $\mathcal{g}_n$ is continuous in the Hausdorff metric at the sphere, then AKC follows. 
%Note that the case for $n=1$ is completed by the equivalence of AKC with (a continuous form of) Conjecture \ref{conj:arithmeticpatterns}, if we consider $1$-dimensional polytopes to be finite subsets of the real numbers.\todo{This is a confusing sentence and I think unnecessary.}

%Just as in the arithmetic Kakeya conjecture, the upper bound of both statements of the conjecture is trivial, and the first statement implies the second.
%For the converse, we slightly modify the corresponding proof for $F_k(N)$ to show that $G_k(N) \leq C \cdot k^2 \big( \log N \big) G_k'(N)$.
%Again, we refer the reader to REF in \cite{GreenRuzsa} for the details.
%Unsurprisingly, \autoref{conj:amcenters} is equivalent to the arithmetic Kakeya conjecture, which we will prove in REF.
%As an intermediate step, we consider \textit{harmonic patterns}.

Finally, in \autoref{sec:fields}, we consider the pattern problem in rings with a more algebraic flavor---finite fields and algebraic number fields.
In \cite{Dvir}, Dvir provides a concise proof of the finite field Kakeya conjecture using the polynomial method.
There is hope that this method can be used to solve Kakeya-type problems by considering their finite field analogs.
For example, the authors in \cite{GreenRuzsa} show that AKC is equivalent to its finite field analog; we show that the same holds for the dual conjecture, \autoref{conj:AMpatterns}.
More generally, we show that many pattern problems in finite fields are (unsurprisingly) equivalent to their analogs in $\Z$.

Green and Ruzsa also use properties of prime numbers to prove the upper bound $F_k'(N) \leq N^{1 - c/\log \log k}$.
We extend this idea to free $\Z$-modules to obtain a similar bound.
\begin{proposition}\label{prop:Gkbound}
For some absolute constant $c$,
\[
G_k'(N) \leq N^{1 - c/\log k}.
\]
\end{proposition}
We also obtain a corresponding result for the plane. Keleti asked in \cite{Squares} and \cite{survey} whether it is possible to reduce the value of $\g(\P)$ if we allow for rotations of $\P$ as well as scaling. In the case of Hausdorff dimension, this problem was studied in \cite{Chang}. Here we show limit result for Minkowski dimension:
\begin{proposition}\label{prop:latticekgon}
    There is a sequence of lattice $k$-gons $\P_k$ that converge in the Hausdorff metric to the unit circle, such that the following holds.
    There exist sets $B_k\subset\R^2$ containing a scaled and rotated copy of the vertices of $\P_k$ centered at every point in $[0,1]^2$ such that
\[ \dim_M B_k \leq 2-\frac{c}{\log k} \]
for some absolute constant $c$.
\end{proposition}
Note that this is a result for \textit{finite} patterns, as $B_k$ need only contain the vertices of a rotated homothet of $\P_k$.

A closing remark: we are indebted to the paper of Green and Ruzsa \cite{GreenRuzsa} as motivation for ours.
Some of our proofs will be almost exactly the same as theirs, only with slight modifications.
In order to keep our paper self-contained, we will typically present these proofs in the appendix.

\subsection{Notation}
For $S \subset \R$, we denote by $S^{+} :=\{ s \in S: s > 0\}$ and follow the convention that $0 \in \N$.
$[n]$ denotes the set $\{1, \dots, n\}$ for $n \in \N$.
Given a function $f$ on $\N$, $f([n])$ denotes the set $\{f(i): i \in [n] \}$.
$U_k$ will typically denote a set of size $k$, and if $U_k \subset R^d$ for some ring $R$, we assume that no coordinate of $U_k$ is $0$.

$\dim_H, \dim_P,$ and $\dim_M$ denote Hausdorff, packing, and Minkowski dimension respectively.
Minkowski dimension is also referred to as \textit{upper box dimension}, and packing dimension as \textit{modified upper box dimension} for the following reason.
For $E \subset \R^n$, let $N_E(\epsilon)$ be the smallest positive integer so that $E$ can be covered by $N_E(\epsilon)$ $n$-dimensional boxes of side-length at most $\epsilon$.
Then,
\[
\dim_M(E) = \limsup_{\epsilon \to 0}\frac{\log N_E(\epsilon)}{\log(1/\epsilon)}.
\]
Further,
\[
\dim_P(E) = \inf \Big\{ \sup_{j} \dim_M(E_j): E \subset \bigcup_{j \in \N}E_j \Big\}.
\]
We will frequently use the following well-known properties.
\begin{enumerate}[(i)]
    \item $\dim_H E \leq \dim_P E \leq \dim_M E$.
    \item If $f: \R^n \to \R^m$ is a Lipschitz map, then $\dim(f(E)) \leq \dim(E)$ for $\dim = \dim_M$ or $\dim_P$.
    \item For $\dim = \dim_M$ or $\dim_P$, $\dim(A \times B) \leq \dim(A) + \dim(B)$ with equality if $B = [0,1]$ or $B = \R$.
    \item $\dim_M\Big( \bigcup_{i=1}^{k}E_i \Big) = \max_{i = 1, \dots, k} \dim_M E_i$, i.e.\ Minkowski dimension is finitely stable.
    \item $\dim_P \Big( \bigcup_{i=1}^{\infty}E_i\Big) = \sup_{i=1, \dots, \infty}\dim_P E_i$, i.e.\ packing dimension is countably stable.
\end{enumerate}
%\todo{This seems more like a number theory paper so it might be helpful to define these things somewhere.}
\section{Finite patterns and arithmetic Kakeya}\label{sec:finitepatterns}

In this section, we consider finite patterns of integers.
\subsection{The pattern conjecture}
%The goal of this subsection is to formulate a family of pattern problems dual to AKC and show that they all imply AKC.
%\begin{definition}
%If $U \subset \Q^{\times}$, a \textit{$U$-pattern with basepoint $x \in \Z$} is a set of the form $x + r(x) \cdot U$, where $r(x)>0$ is called the \textit{scaling factor} of the pattern.
%\end{definition}
%For example, when $U = [k]$, a $U$-pattern is a $k$-term arithmetic progression.
%We are interested in questions of the following type: if $B \subset \R$ contains a $U$-pattern for every basepoint in a set $S$, what can we say about the size of $B$?
Let $(U_k)_{k \in \N}$ be an increasing chain of subsets of $\Q^{+}$ such that $|U_k| = k$.
The following conjecture is technically a family of conjectures, one for each such chain.
\begin{conjecture}[The $(U_k)$-pattern conjecture]\label{conj:pattern}
Let $k$ and $N$ be positive integers.
Let $G'_{U_k}(N)$ be the size of the smallest set of integers that contains a $U_k$-pattern with basepoint $x$, for $N$ distinct integer values of $x$. Then,
\[
\lim_{k\to\infty} \lim_{N\to\infty} \frac{\log G'_{U_k}(N)}{\log N} = 1.
\]
\end{conjecture}
A few remarks are in order.
First, this limit is well-defined as $U_k \subset U_{k+1}$, so $G'_{U_{k+1}}(N) \geq G'_{U_k}(N)$.
Second, we may replace ``$N$ distinct integer values of $x$'' with ``$N$ distinct rational values of $x$'' and obtain the same quantity.
Any set $B \subset \Q$ containing a $U_k$-pattern with $N$ distinct basepoints can be multiplied with a large enough common denominator $c$ so that $c \cdot B \subset \Z$ also contains a $U_k$-pattern with $N$ distinct \textit{integer} basepoints.

In line with AKC, we can formulate a related conjecture.
Now, let $(U_k)_{k \in \N}$ be an increasing chain of subsets of $\N^+$ such that $|U_k| = k$.
\begin{conjecture}[The restricted $(U_k)$-pattern conjecture]
 Let $k$ and $N$ be positive integers.
Let $G_{U_k}(N)$ be the size of the smallest set of integers that contains a $U_k$-pattern with basepoint $x$, for each $x \in [N]$. 
Then,
\[
\lim_{k\to\infty} \lim_{N\to\infty} \frac{\log G_{U_k}(N)}{\log N} = 1.
\]
\end{conjecture}
As expected, for sequences that satisfy $U_k \subset \N^{+}$, the pattern conjecture and the restricted pattern conjecture are equivalent.
\begin{proposition}\label{prop:restrictedequiv}
Fix an increasing chain $(U_k)_{k \in \N^+}$ of subsets of $\N$ such that $|U_k|=k$.
Then,
\[
\lim_{k\to\infty} \lim_{N\to\infty} \frac{\log G_{U_k}(N)}{\log N} = \lim_{k\to\infty} \lim_{N\to\infty} \frac{\log G'_{U_k}(N)}{\log N}.
\]
\end{proposition}
We relegate the proof to the appendix as it is almost exactly the same as one presented by Green and Ruzsa in \cite{GreenRuzsa}.
When $U_k = [k]$, we recover the equivalence in the arithmetic pattern conjecture, as stated in \autoref{conj:AMpatterns}.
For convenience, we will write ``the pattern (resp.\ restricted pattern) conjecture'' instead of ``for each sequence $(U_k)_{k \in \N}$, the $(U_k)$-pattern (resp.\ restricted $(U_k)$-pattern) conjecture''.

We will now show that the pattern conjecture implies AKC. 
%By \autoref{prop:restrictedequiv}, the restricted pattern conjecture also implies AKC.
For the proof, we will need another equivalent formulation of AKC.
For $A \subset \Z \times \Z$ and $r \in \Q \cup \{\infty\}$, define
\[
\pi_r(A) = \{ x + r y: (x,y) \in A \},
\]
with the convention that $x + \infty y = y$.
\begin{conjecture}[AKC projection formulation, \cite{GreenRuzsa}]\label{conj:projection}
For $\epsilon > 0$ arbitrary, there are $r_1, \dots, r_k\in \Q \cup \{\infty\},$ none equal to $-1$, such that for any $A \subset \Z \times \Z$,
\[
| \pi_{-1}(A) | \leq \max_{i=1, \dots, k} | \pi_{r_i}(A) |^{1+\epsilon}.
\]
\end{conjecture}
The main tool in our proof will be a modification of a lemma from \cite{GreenRuzsa}, the proof of which we again relegate to the appendix.

\begin{lemma}[Amplification lemma]
Suppose that there exist $B\subset \Z \times \Z$, a set $S \subset \N$ of size $k$, and $\epsilon>0$ such that
\[
\lvert \pi_{-1}(B) \rvert > \max_{j \in S} \lvert \pi_j(B) \rvert^{1 + \epsilon}.
\]
Then, for all $M>0$, there exists $B'\subset\Z \times \Z\setminus\{0\}$ such that
\[
\lvert \pi_{-1}(B') \rvert > M \max_{j\in S} \lvert \pi_j(B') \rvert^{1 + \epsilon}.
\]
\end{lemma}

\begin{theorem}\label{thm:patternimpliesakc}
The pattern conjecture implies AKC.
\end{theorem}
\begin{proof}
Fix a sequence $(U_k)_{k \in \N}$ and suppose that the projection formulation of AKC fails. 
Then, there exists $\epsilon > 0$ such that for any $r_1, \hdots r_k\in \Q \cup\{\infty\}\setminus\{-1\}$ there is a set $B \subset \Z \times \Z$ so that
\[
| \pi_{-1}(B)| > \max_{i=1, \dots, k} |\pi_{r_i}(B)|^{1+\epsilon}.
\]
Let $U_k - 1 = \{ u-1: u \in U_k\}$.
Since $0 \notin U_k$, $-1 \notin U_k - 1$.
For each $k$, let $m > 0$ be arbitrary, and define $M_k=mk^{1+\epsilon}$. 
Then, applying the amplification lemma for $B$ and $M_k$, there exist sets $B_{k,m}\subset \Z \times \Z$ such that
\[
\lvert \pi_{-1}(B_{k,m}) \rvert >m\left(k\max_{r\in U_k-1}\lvert \pi_r(B_{k,m})\rvert \right)^{1+\epsilon}.
\]
By translating $B_{k,m}$ if necessary, we may assume $B_{k,m}\in\Z\times\Z^+$ while maintaining the size of its projections. Define
\[
A_{k,m} := \bigcup_{r\in U_k-1} \pi_r(B_{k,m}) \subset \Q.
\]
Since $x+y\cdot(U_k-1)=(x-y)+y\cdot U_k$, $A_{k,m}$ contains a $U_k$-pattern with basepoint $x-y$ and scaling factor $y$ (recall that $y$ is nonzero by construction) for every $(x,y) \in B_{k,m}$.

On the other hand,
\[
k \leq \lvert A_{k,m} \rvert \leq k \max_{j \in S} \lvert  \pi_j(B_{k,m}) \rvert < \Big( \lvert \pi_{-1}(B_{k,m}) \rvert/m \Big)^{1/(1+ \epsilon)}.
\]
Let $N_m = \lvert \pi_{-1}(B_{k,m}) \rvert$, so the above inequalities yield $N_m > mk^{1+\epsilon}$.
Since $m$ was arbitrary, $N_m \xrightarrow[]{m\to\infty}\infty$, so
\[
\lim_{N_m \to \infty} \frac{\log G_{U_k}'(N_m)}{\log N_m} \leq \lim_{m \to \infty}\frac{\log N_m^{1/(1+\epsilon)}}{\log N_m} = \frac{1}{1+ \epsilon} < 1.
\]
Since this is a uniform bound for all $k$, the pattern conjecture fails.
\end{proof}

\subsection{The arithmetic and harmonic pattern conjectures}

Perhaps the most interesting cases of the pattern conjecture are when $U_k = [k]$ (the arithmetic pattern conjecture) and when $U_k = 1/[k]$ (the harmonic pattern conjecture).
In this subsection, we prove \autoref{thm:3patterns}.
Recall that $1/[k] = \{ 1/i: i \in [k] \}$.

\begin{theorem}\label{thm:harmonicpattern}
For $U_k = 1/[k]$, the pattern conjecture is equivalent to the Arithmetic Kakeya conjecture.
In fact,
\[
\frac{1}{k}\cdot G'_{1/[k]}(N) \leq F_k'(N) \leq k \cdot G'_{1/[k]}(N).
\]
\end{theorem}
\begin{proof}
Suppose $A$ contains a $1/[k]$-pattern for every basepoint $x \in S$, with corresponding scaling factors $r(x)$.
Define
\[
A_i = \Big\{ x + r(x) \cdot \frac{1}{i}: x \in S \Big\}.
\]
Then, the set $B = \bigcup_{i=1}^{k}i \cdot A_i$ contains a $k$-term arithmetic progression for every common difference $x \in S$, and $|B| \leq k \cdot |A|$.
So, $F_k'(N) \leq k \cdot G'_{1/[k]}(N)$.

Similarly, if $B$ contains a $k$-term arithmetic progression for every common difference in $S$, let $B_i \subset B$ be the set of the $i$th terms of the progressions.
The set $A = \bigcup_{i}1/i \cdot B_i$ contains a $1/[k]$-pattern for every basepoint in $S$, and $|A| \leq k \cdot |B|$, completing the proof.
\end{proof}

Now,
\begin{theorem}\label{thm:AMPimpliesAKC}
The arithmetic pattern conjecture is equivalent to AKC.
\end{theorem}
\begin{proof}
By \autoref{thm:patternimpliesakc} and \autoref{thm:harmonicpattern}, we only need to show that the $1/[k]$-pattern conjecture implies the $[k]$-pattern conjecture.
Suppose $B \subset \Z$ contains a $k!$-term progression for every basepoint in $S \subset \Z$, i.e.
\[
B = \bigcup_{x \in S} x + r(x) \cdot [k!].
\]
We claim that $B$ contains a $1/[k]$-pattern for every basepoint in $(c/k!)\cdot S$.
That is,
\[
\bigcup_{x \in S} x + (k! \cdot r(x)) \cdot 1/[k] \subset B.
\]
So,
\[
G_{[k!]}'(N) \geq G_{1/[k]}'(N) \geq \frac{1}{k} \cdot F_k'(N),
\]
and this completes the proof.
\end{proof}
The key tool in this proof was observing that the set $1/[k]$ can be multiplied into the set $[k!]$: every $[k!]$-pattern contains a $(1/[k])$-pattern.
%By a similar manipulation, we can prove the converse implication of \autoref{thm:patternimpliesakc} for a larger family of sequences $(U_k)_{k \in \N}$.
%\begin{proposition} AKC is equivalent to the $(U_k)$-pattern conjecture for the following choices of $(U_k)$. \begin{enumerate}[(a)]  \item Let $f:\N \to \N$ be a bijection. Set $U_k = f([k])$.   \item Let $(V_k)_{k \in \N}$ be subsets of $\N$ of cardinality $k$ which are bounded in value by $O(k)$. Then there exist integers $t_k \in \mathbb{Z}$ so that we may take $U_k = V_k + t_k$.\end{enumerate}\end{proposition}
%\begin{proof}
%For (a), set $n_k = \max f^{-1}([k])$.
%Then, $[k] \subset f([n_k])$, so $G'_{U_{n_k}}(N) \geq G'_{k}(N)$.

%The proof of (b) uses the finitary version of Szemer\'edi's theorem for arithmetic progressions.
%Namely, fix $k \in \N$.
%If $V_n \subset [cn]$ for some absolute constant $c$, then for some $n_k$ sufficiently large, $V_{n_k}$ contains a $k$-term arithmetic progression, say $t_{n_k} + r(t_{n_k}) \cdot [k]$.
%Setting $U_{n_k} = V_{n_k} - t_k$, we mimic the proof of (a) to show that $G'_{U_{n_k}}(N) \geq G'_{k}(N)$.
%\end{proof}

\subsection{Bounds and constructions}
In this subsection, we consider nontrivial estimates for some patterns.
We will provide upper bounds in later sections, once we have developed the necessary tools.
In \autoref{sec:fields}, we use a probabilistic construction to show that $G_{U_k}'(N) \leq N^{1-1/k + o_N(1)}$ for any sequence $(U_k)_{k \in \N}$.
If we assume that the sequence $(U_k)_{k \in \N}$ grows subexponentially in $k$, we use algebraic number theory to improve this to $G_{U_k}'(N) \leq N^{1- c/\log k}$.

However, these bounds are trivial for small values of $k$, so we study the behavior of $G_k'(N)$ when $k$ is small.
A trivial sumset estimate yields $G_2'(N) = \Theta(N^{1/2})$, and this is sharp when 
\[
B = \{ 2^n: 0 \leq n \leq N-1 \},
\] 
and 
\[
S = \{x \in \Z: B \text{ contains a }2\text{-term arithmetic progression centered at }x\}.
\]

\begin{example}
For any integers $N$ and $k$, $G_k'(N) =O( N \log N)$.
This will be a sort of reverse-engineered construction.
Fix $k$, and for every $n$, we construct a set $B_{n}$ of size $k^n$ so that $B_n$ contains a $k$-term progression for every basepoint in a set $S_n$ of size $nk^{n-1}$.
Let $B_1 = [k]$ and $S_1 = \{ 0 \}$.
If $B_n$ and $S_n$ have been constructed according to our hypothesis, choose $t \in \N$ large enough so that
\[
\min(t + S_n) > \max(B_n).
\]
Now, set
\[
B_{n+1} = \bigcup_{j=0}^{k-1}\Big(j \cdot t + B_n\Big), \quad
S_{n+1} = \bigcup_{j=0}^{k-1}\Big(j \cdot t + S_n\Big) \bigcup \Big( -t + B_n \Big).
\]
It is routine to check that every point of $B_{n+1}$ contains a $k$-term progression for every basepoint in $S_{n+1}$. 
By our choice of $t$, $|B_{n+1}| = k^{n+1}$ and $|S_{n+1}| = (n+1)k^n$.
\end{example}

Now we consider lower bounds on $G_k'(N)$.
Combining \autoref{thm:AMPimpliesAKC} with the bound $F_3'(N) \geq N^{6/11}$ by Katz and Tao (\cite{KatzTao}), we obtain
\begin{proposition}
\[
G_{6}'(N) \geq \frac{1}{3}N^{6/11}.
\]
\end{proposition}
We can get rid of the cumbersome factor of $1/3$ by directly modifying the proof of Katz and Tao.
This, in fact, proves something slightly stronger.
\begin{proposition}\label{prop:346}
If $U = \{3,4,6\}$, then $G_{U}'(N) \geq N^{6/11}$.
\end{proposition}
For this, we apply a result of Katz and Tao.
Let $A, B \subset \Z$, $G \subset A \times B$, and denote by $A \stackrel{G}{+} B = \{ a+b: (a,b) \in G \}$.
\begin{theorem}[\cite{KatzTao}]\label{thm:katztao}
Suppose $A_1, A_2 \subset \Z$, $G \subset A_1 \times A_2$, and $|A_1|, |A_2|, |A_1 \stackrel{G}{+} A_2| \leq n$.
Then, $|A_1 \stackrel{G}{-}  A_2| \leq n^{11/6}$.
\end{theorem}
\begin{proof}[Proof of \autoref{prop:346}]
Suppose $B = \bigcup_{x \in S} x + r(x) \cdot U$, where $S \subset \N$ is a set of size $N$.
Define
\[
A_1 = \Big\{ x + 6r(x): x \in S \Big\},\quad A_2 = \Big\{ x + 3r(x): x \in S \Big\},
\]
and
\[
G = \Big\{ \big( x+6r(x), 2x + 6r(x) \big): x \in S \Big\} \subset A_1 \times 2A_2.
\]
The set $A_1 \stackrel{G}{+} 2A_2$ is in bijection with $\{x + 4r(x): x \in S\}$, so applying \autoref{thm:katztao} yields $|B| \geq |S|^{6/11}$.
\end{proof}

In general, we can use \autoref{thm:katztao} to find nontrivial bounds for $G'_U(N)$, when $U$ has some structure. 
%Our main result from this investigation is the following bound.
%\begin{corollary}
%\[
%G_4'(N) \geq N^{4/7}.
%\]
%\end{corollary}
For example, we consider $3$-term patterns.

\begin{proposition}\label{prop:3term}
If $U$ has the form $U = \{pq,q,2pq/(p+1) \}$, then $G_{U}'(N) \geq N^{6/11}$.
%For any set $P = \{ d_1, d_2, d_3 \} \subset \Q \setminus \{0\}$ with $d_1 < d_2 < d_3$, there exists some $c \in \Q$ such that if $B$ contains a $(P-c)$-pattern for every basepoint in a set $S \subset \Z$, then $|B| \geq |S|^{6/11}$.
\end{proposition}
\begin{proof}
%The proof is a rather technical manipulation.
%Set
%\begin{align*}
  %  c &= d_1 + \frac{(d_3 - d_1)(d_2 - d_1)}{(2d_3 - d_1 - d_2)}, \\
   % p &= - \frac{d_3 - d_1}{d_3 - d_2}, \\
  %  q &= \frac{(d_2 - d_1)(d_3 - d_2)}{2d_3 - d_1 - d_2}.
%\end{align*}
%For example, if $P = \{1,2,3\}$, then $P-c = \{-2/3, 1/3, 4/3\}$.
%In general, for $P =\{d_1, d_2, d_3\}$, we get $P - c = \{pq, q, 2pq/(p+1)\}$. 
Suppose $B$ contains a $U$-pattern for every basepoint in a set $S \subset \Z$.
Define 
\[
A_1 = \{ x+ r(x)pq: x \in S \}, \quad A_2 = \{ x + r(x)q: x \in S \}
\]
and
\[
G = \Big\{ (x + r(x)pq, x + r(x)q) : x \in S \Big\} \subset A_1 \times A_2.
\]
Then,
\[
A_1 \stackrel{G}{+} p A_2= (p+1)\Big\{ x+ r(x)\frac{2pq}{p+1}: x \in S \Big\}, \quad A_1 \stackrel{G}{-} p A_2 = (1-p)S.
\]
Notice, that $A_1, A_2,$ and $(1/p+1) \cdot \big( A_1 \stackrel{G}{+} pA_2 \big) \subset B$, so $|A_1|, |A_2|, |A_1 \stackrel{G}{+} p A_2| \le |B|$. 
Applying \autoref{thm:katztao}, $|B| \geq |S|^{6/11}$.
\end{proof}

\section{Patterns in the continuous setting}\label{sec:infinitepatterns}
\subsection{Dimension considerations}
In this section we consider pattern problems where the correct notion of size is dimension.
Our first main result is that for problems of this form, Minkowski and packing dimensions are essentially equivalent, while Hausdorff dimension gives little information.
Recalling our notation from the introduction, let
\[
\cal{P} = \{ \P_{x,r}: x \in X, r \in R \} 
\]
be a double-indexed family of subsets of $\R^n$, such that each $\P_{x,r}$ is finite.
As usual, we will typically be interested in the case $\P_{x,r} = x + r \cdot \P$, for some set $\P \subset \R$.

For $B \subset \R^n$, define
\[
S(B, \cal{P}):= \{ r \in R : \P_{x,r} \subset B \text{ for some } x \in X \}.
\]
Informally, this is the set of scaling factors for which $B$ contains a $\P$-pattern.
With this notation, we can talk about the equivalence of Minkowski and packing dimensions in the following sense.
\begin{theorem}\label{thm:box=packing}
Suppose there is an absolute constant $c > 0$ such that for any bounded subset $B \subset \R^n$, $\dim_M B \geq c \cdot \dim_MS(B, \cal{P})$.
Then, for any subset $B' \subset \R^n$, $\dim_P B' \geq c \cdot \dim_P S(B', \cal{P})$.
\end{theorem}

\begin{corollary}
Fix a finite set $\P \subset \R$.
Let $\f_n^{M}(\P)$ and $\f_n^{P}(\P)$ denote the minimal Minkowski and packing dimensions respectively of a set $B$ containing a $\P$-pattern for every scaling factor in $[0,1]^n$.
Then, $\f^{M}_n(\P) = \f^{P}_n(\P)$.
\end{corollary}

\begin{corollary}
Fix a finite set $\P \subset \R$.
Let $\g_n^{M}(\P)$ and $\g_n^{P}(\P)$ denote the minimal Minkowski and packing dimensions respectively of a set $B$ containing a $\P$-pattern for every basepoint in $[0,1]^n$.
Then, $\g_n^{M}(\P) = \g_n^{P}(\P)$.
\end{corollary}

\begin{proof}[Proof of \autoref{thm:box=packing}]
Fix a set $B\subset \R^n$, and let $\{B_j\}_{j=1}^\infty$ be any collection of subsets of $\R^n$ such that $B = \bigcup_{j} B_j$.
Denote $B'_m=\cup_{j=1}^{m} B_j$. 
By assumption, we know that 
    \[
    \dim_M B'_m \ge c \cdot \dim_M S(B'_m,\cal{P}).
    \]
    The sets $B_m'$ form an ascending chain of sets such that $\bigcup_{m}B_m' = B$.
    If $\P_{x,r} \subset B$, since it is finite, $\P_{x,r} \subset B'_m$ for some large enough $m$, i.e.\ $S(B,\cal{P})=\bigcup_{m}S(B'_m,\cal{P})$.
    By  the definition of packing dimension, $\dim _ P S(B, \cal{P}) \leq \sup_{m} \dim_M S(B_m', \cal{P})$.
    Hence, for every $\epsilon>0$, there is some $m$ such that
\[
\dim_M S(B'_m,\cal{P})>\dim_P S(B,\cal{P})-\epsilon.
\]
Then,
\[
\max_{j=1, \dots, m} \dim_M B_j=\dim_M B'_m\geq c(\dim_P S(B,\cal{P})-\epsilon).
\]
Since $\epsilon$ was arbitrary, $\sup_{j \in \N} \dim_M B_j \geq c \cdot \dim_P S(B,\cal{P})$, and thus $\dim_P B \ge c \cdot \dim_P S(B,\cal{P})$.
\end{proof}

On the other hand, Hausdorff dimension fails to provide much information due to the following result. 

%Unfortunately, Hausdorff dimension does not tell us anything.
\begin{proposition}[Proposition 3.1, \cite{Squares}]\label{prop:hausdorff}
For any dimension $n \geq 1$, there is a closed set $A \subset \R^n$ with $\dim_H A = 0$, such that for any finite family of invertible affine maps $(f_i)_{i=1}^{k}$ on $\R^n$, the intersection $\bigcap_{i=1}^{k}f_i(A)$ is nonempty.
If the maps $f_i$ are constrained to lie in a fixed compact set of invertible affine maps, $A$ can be chosen to be compact.
\end{proposition}
For example, for $i \in \N$ and $r \in [0,1]^n$, define $f_{i,r}(a) = a + i \cdot r$.
Then there is a set $A$ containing an infinite arithmetic progression for every common difference $r \in [0,1]^n$, with $\dim_H A = 0$.
If we restrict to $k$-term arithmetic progressions, so $i \in [k]$, then $A$ can be chosen to be compact.

\subsection{Higher dimensional pattern conjectures}
Here, we consider three generalizations of AKC.
The astute reader will now realize that they must be equivalent.
Each of the following is technically a family of conjectures, one for each positive integer $n$.
\begin{conjecture}[$n$-dimensional AKC]\label{conj:ndimAKC}
Let $F'_{n,k}(N)$ denote the minimal cardinality of a set containing a $k$-term arithmetic progression in $\mathbb{R}^n$ with $N^n$ distinct common differences in $\N^n$. 
Then
\[
\lim_{N \to \infty} \lim_{k \to \infty} \frac{\log F'_{n,k} (N)}{\log N} = n.
\]
\end{conjecture}
For convenience, we will denote the limit $\lim_{N \to \infty}\big(\log F'_{n,k} (N)/ \log N \big)$ as $F'_{n,k}$.
\begin{conjecture}[$n$-dimensional Minkowski AKC]\label{conj:boxAKC}
 Let $\mathcal{f}^M_n (k)$ be the infimum of the Minkowski dimension of a set $B_k$ in $\R^n$ that contains a $k$-term arithmetic progression with every common difference in $[0,1]^n$.
Then,
\[
\lim_{k \to \infty}\mathcal{f}^M_n (k) = n.
\]
\end{conjecture}
\begin{conjecture}[$n$-dimensional packing AKC]\label{conj:packingAKC}
 Let $\mathcal{f}^P_n (k)$ be the infimum of the packing dimension of a set $B_k$ in $\R^n$ that contains a $k$-term arithmetic progression with every common difference in $[0,1]^n$.
Then, 
\[
\lim_{k \to \infty}\mathcal{f}^P_n (k) = n.
\]
\end{conjecture}
When $n=1$, we simply refer to the corresponding statements as the Minkowski and packing dimension arithmetic Kakeya conjectures. 
Our main result for this subsection is that
\begin{theorem}\label{thm:AKCequivs}
The above conjectures are equivalent for all $n$. 
Moreover,
\[
\f_n^{P}(k) = \f_n^{M}(k) = F'_{n,k},
\]
that is, the conjectures in the three \textit{families} are all equivalent.
\end{theorem}
Of course, an immediate corollary of \autoref{thm:box=packing} is that $\f_n^{M}(k) = \f_n^{P}(k)$. We remark that our theorem is stronger than stating that ``in each dimension, the three conjectures are equivalent'': the conjectures in the three \textit{families} are all equivalent.
%\begin{conjecture}[Minkowski dimension arithmetic Kakeya]\label{conj:ctsamkakeya-box}
%Let $\mathcal{f}_B(k)$ be the infimum of the Minkowski dimension of a set $B_k$ in $\R$ that contains a $k$-term arithmetic progression with every common difference in $[0,1]$.
%Then,
%\[
%\lim_{k \to \infty}\mathcal{f}_B(k) = 1.
%\]
%\end{conjecture}
%\begin{conjecture}[Packing dimension arithmetic Kakeya]\label{conj:ctsamkakeya-packing}
 %Let $\mathcal{f}_P(k)$ be the infimum of the packing dimension of a set $B_k$ in $\R$ that contains a $k$-term arithmetic progression with every common difference in $[0,1]$.
%Then,
%\[
%\lim_{k \to \infty}\mathcal{f}_P(k) = 1.
%\]
%\end{conjecture}

The crucial step is a modification of a proof of \cite{GreenRuzsa} (as usual, to be found in the appendix).
\begin{theorem}\label{thm:ndimAKC}
     %\autoref{conj:AKC} is equivalent to \autoref{conj:ndimAKC} for each $n>0$.
     The Arithmetic Kakeya conjecture is equivalent to the $n$-dimensional AKC for each $n > 0$.
\end{theorem}

%For the rest of this section, our goal is to prove the equivalence of these three conjecture. 
%\begin{proof}
%Set $n=1$, and let $A$ be a set as in \autoref{prop:hausdorff}.
%Let $f_{i,d}(x) = x - i \cdot d$ these are invertible affine maps.
%Then, for any $k \in \N$ and $d \in [0,1]$, $\bigcap_{i=1}^{k}f_{i,d}(A) \neq \emptyset$, so $A$ contains a $k$-term arithmetic progression with common difference $d$.
%If we only require this for finitely many $k$, the maps $f_{i,d}: d \in [0,1]$ lie in a compact set, and the second statement follows.
%\end{proof}
%For convenience, we will denote the limit $\lim_{N \to \infty}\big(\log F_{n,k} (N)/ \log N \big)$ as $F_{n,k}$.
%Now, we will state Hutchinson's theorem (Theorem 9.3 of \cite{Falconer}) and two other results used in the proof of the equivalences. 
Next, we will need some preliminary results.
\begin{lemma}[\cite{Falconer}]\label{lem:falconer}
Let $\phi_1, \dots, \phi_m$ be similarities on $D \subset \R^n$ with ratios $c_1, 
\dots, c_m$ respectively such that $c_i < 1$ for each $i$.
For any $F \subset \R^n$, define $\phi(F) = \bigcup_{i=1}^{m}\phi_i(F)$.
Suppose $\phi$ satisfies the \textit{open set condition}: for some nonempty open set $V$, $\phi(V) \subset V$.
Then, if $F$ is a nonempty compact set such that $F = \phi(F)$, then $\dim_H F = \dim_M F = s$, where $\sum_{i=1}^{m}c_i^s = 1$.
\end{lemma}
\begin{lemma}\label{lem:boxAKC} For any $n$,
\[
\mathcal{f}_n^M(k)\geq  F'_{n,k}.
\]
\end{lemma}
\begin{proof}
Let $P$ be a set containing a $k$- term arithmetic progression for every common difference in $[0,1]^n$ such that $\dim_M P< \f_n^{M}(k) + \epsilon$. 
We can write
\[
P = \Big\{ a(r) + i \cdot r: r \in [0,1]^{n}, 0 \leq i \leq k-1 \Big\}.
\]
For each $q \in \N$, subdivide $\R$ into disjoint half-open intervals of length $1/q$.
For $r=(r_1, \hdots,r_n)\in [0,1]^n$, denote the right endpoint of the interval containing $r_j$ as $r_{j}(q)$, and set $r(q)=(r_{1}(q), \hdots,r_{n}(q))$. 
Notice that $\Big| \{r(q) : r\in [0,1]^n\} \Big| = q^{(1+o(1))n}$.
Define $a(r;q)_j$ and $a(r;q)$ analogously.
This map sends each point to the right corner of the $(1/q)$-cube containing it.

Set 
\[
E_{q,i} = \Big\{ a(r;q) + i\cdot r(q): r \in S \Big\} 
\]
for each $0 \leq i \leq k-1$.
The $(1/q)$-cubes determined by the points of $E_{q,i}$ are disjoint sets containing points of $P$, so
\[
\frac{\log |E_{q,i}|}{\log q} \leq (1+o_q(1))\dim_M P \leq (1+o_q(1))\f_n^M(k)+ \epsilon.
\]
%Hence, $|E_{q,i}|\leq q^{(1+o_q(1))\f_n^M(k)+ \epsilon}$.
Define 
\[
E_q=\{a(r;q)+i \cdot r(q): r\in [0,1]^n, 0 \leq i \leq k-1\}; \quad |E_q| \leq k\max_{0 \leq i \leq k-1}|E_{q,i}|.
\]
Setting $S_q=q \cdot E_q$, we have that $S_q \subset \Z^n$ contains a $k$-term arithmetic progression for $N^n$ distinct differences in $\Z^n$, and 
 \[
 \lim_{q \to \infty} \frac{\log |S_q|}{\log q} \leq  \mathcal{f}_n^M (k)+ \epsilon.
\]
Since this holds for any $\epsilon>0$, 
\[
\lim_{q\to\infty}\frac{\log F'_{n,k} (q)}{\log q}\leq \mathcal{f}_n^M(k).
\]
\end{proof}
Now we prove \autoref{thm:AKCequivs}.
\begin{proof}[Proof of \autoref{thm:AKCequivs}]
By \autoref{thm:ndimAKC}, it suffices to prove the equivalence in a fixed dimension $n$.
We already know that $\f_n^{M}(k) = \f^{P}_n(k)$, and that $\f^M_n(k) \geq F'_{n,k}$.
It remains to show that $\f^M_n(k) \leq F'_{n,k}$. 
Let $A \subset \Z^n$ be a set of size $F'_{n,k}(N)$ containing a $k$-term arithmetic progression for $N^n$ distinct common differences.
Define
\[
F = \Big\{\sum_{j=0}^{\infty} a_jN^{-j}: a_j \in A \Big\}.
\]
For each $a \in A$, define the map $\phi_a : \R^n \to \R^n$ by $\phi_a(x) = a + (x/N)$; this is a similarity with ratio $1/N$.
It is easy to check that $F$ satisfies the conditions of \autoref{lem:falconer}, so $\dim_M F \leq \log F'_{n,k}(N)/\log N$.
Further, $F$ contains a $k$-term arithmetic progression for every common difference $r \in [0,1]^n$. 
Namely, for $r=(r_1, r_2, \hdots, r_n)\in [0,1]^n$,
write each $r_j$ in base $N$ as $r_j = \sum_{m=0}^{\infty} r_{j,m} N^{-m}$.

Then, the points
\[
\sum_{m=0}^{\infty}
\begin{bmatrix}
a(r_{1,m})N^{-m} \\
a(r_{2,m})N^{-m} \\
\hdots\\
a(r_{n,m})N^{-m} 
\end{bmatrix}
+ i \cdot \sum_{m=0}^{\infty}
\begin{bmatrix}
r_{1,m}N^{-m}\\
r_{2,m}N^{-m}\\
\hdots\\
r_{n,m}N^{-m}
\end{bmatrix}
: 1 \leq i \leq k
\]
lie in $F$ and form the desired arithmetic progression.
Since this holds for every $N$,
\[
\f_n^{M}(k) \leq F'_{n,k}
\]
and this concludes the proof.
\end{proof}

Specifically, when $n=1$, we recover \autoref{thm:1dAKCequivs} that $\f_1^{M}(k) = \f_1^{P}(k) = F'_{1,k}$. 
The proof of \autoref{thm:AKCequivs} can be replicated for continuous one-dimensional analogs of the pattern conjecture.

\begin{theorem}
Define $\g^M(U_k)$ (resp.\ $\g^P(U_k)$ to be the minimal Minkowski (resp.\ packing) dimension of a set containing a $U_k$-pattern for every basepoint in $[0,1]$. 
Then,
\[
\g^{M}(U_k) = \g^P(U_k) = \lim_{N \to \infty} \frac{\log G'_{U_k}(N)}{\log N}.
\]
\end{theorem}
%We will use this identity in \autoref{sec:fields} to obtain an upper bound for $G_{U_k}'(N)$.
In \cite{BourgainKakeya}, Bourgain showed that the Arithmetic Kakeya conjecture implies the Kakeya conjecture for Minkowski dimension.
We will show that it in fact implies the Kakeya conjecture for packing dimension.
More specifically,
\begin{proposition}
The $n$-dimensional AKC for packing dimension (\autoref{conj:packingAKC}) implies the Kakeya conjecture for packing dimension in $\R^n$.
\end{proposition}
\begin{proof}
Fix $\epsilon>0$. 
Assuming \autoref{conj:packingAKC} holds, choose $k \in \N$ such that $\mathcal{f}^P_n(k) \geq n-\epsilon$.
Suppose $B \subset \R^n$ contains a line segment of length $k\sqrt{n}$ in each direction.
($B$ is a scaled Besicovitch set.)
Namely, for each unit vector $\hat{u}\in\R^n$, there is $x(\hat{u}) \in \R^n$ such that $B$ contains the line segment $x(\hat{u}) + \hat{u}\cdot [0,k\sqrt{n}]$.
For $u \in [0,1]^n$, let $\hat{u} = u/||u||$.
Then $B$ contains a $k$-term arithmetic progression with every common difference $u\in[0,1]^n$, so $\dim_P B \geq n - \epsilon$.
Since $\epsilon$ was arbitrary, the Kakeya conjecture follows.
\end{proof}

\section{Geometry}\label{sec:geometry}
In this section, we study homothets of convex polytope with centers in a given set $S \subset \R^n$. 
(Recall that by polytope, we mean the boundary of a finite intersection of half-spaces).
We will prove \autoref{thm:harmonicpolytope}, and the bounds listed in \autoref{table:geom}.
%In fact, we will prove something slightly stronger than \autoref{thm:harmonicpolytope}.
Fix a dimension $n>1$.
We fix the following notation.

Suppose $\P$ is a convex polytope in $\R^n$, defined by a collection of bounding hyperplanes, $\{H_j : j\in J\}$ for some index set $J$. 
We say $\P$ is a $m$-tope if it has $m$ bounding hyperplanes.
The $j$\textsuperscript{th} hyperplane $H_j$ can be defined by the pair $(\vec{u}_j, d_j)$, where $\vec{u}_j$ is the unit normal vector to $H_j$, and $d_j$ the distance of $H_j$ from the origin.
If $x + r \cdot \P$ is a homothet of $\P$ centered at $x$, we denote the corresponding translate of the hyperplane $H_j$ by $H_j(x)$.

\subsection{Polytope sequences and AKC}

The following theorem is technically a family of theorems for each dimension.
Recall that $\f(m) = \f^{P}(m) = \f^{M}(m)$ is the minimal packing or Minkowski dimension of a set in $\R$ that contains an $m$-term arithmetic progression for every common difference in $[0,1]$.
These infima are equal for packing and Minkowski dimensions by \autoref{thm:box=packing}.
\begin{theorem}[$n$]\label{thm:harmonicpolytopestrong}
With the above notation, let $\P\subset\R^n$ be an polytope with the property that $d_j = (\vec{u}_j)_1/j$ for all $j\in J$, where the index set $J$ is an arithmetic progression of length $m$. Let $\mathcal{h}^M(\P)$ and $\mathcal{h}^P(\P)$ denote the minimal Minkowski and packing dimensions (respectively) of a set containing a homothet of $\P$ centered at each point in the straight line segment $[0,1]\times\{0\}^{n-1}$. Then 
\[\mathcal{h}_n^M(\P) = \mathcal{h}_n^P(\P) = n-1+\mathcal{f}(m).\]
\end{theorem}
\begin{proof}[Proof of Theorem \ref{thm:harmonicpolytopestrong}]

Denote $(x, 0 ,\dots, 0) \in [0,1] \times \{0\}^{n-1}$ by $\hat{x}$. We first show that $\h_n^P(\P)\geq n-1 + \f^P(m)$. Suppose $B$ is a set containing a homothetic copy of $\P$ around each $\hat{x}\in[0,1]\times \{0\}^{n-1}$, so that for each $x\in[0,1]$, $\hat{x} + r(x) \cdot \P\in B$.
Define
\[
A_{j} = \Big\{ \langle \vec{x}, \vec{u}_j \rangle + d_j \cdot r(x) : x \in [0,1] \Big\} = \Big\{ x \cdot (\vec{u}_j)_1  + d_j \cdot r(x) : x \in [0,1] \Big\} \subset \R.
\]
An elementary calculation shows that this quantity is the distance of $H_j(\hat{x})$ from the origin.
We now use the following standard trick.
Let $F_j$ be the $(n-1)$-face of $
\P$ contained in $H_j$, and let $v_1, \dots, v_{n-1} \in \R^n$ be a basis for the linear subspace $U_j$ corresponding to $H_j$.
It follows that
\[
H_j = F_j + \sum_{i=1}^{n-1}\Q \cdot v_{i}.
\]
Repeating this argument for each hyperplane, it follows that there is a countable set $C$ such that $B' = B+C$ contains the entire hyperplane $H_j(\hat{x})$ for each $x \in [0,1]^{n}$.
%Each polytope is a union of faces, which are pieces of hyperplanes. In the next step, we show that packing dimension is invariant if we replace these hyperplane pieces by entire hyperplanes. Each face $F_j$ of $\P$ is an $n-1$ dimensional polytope, so in particular, it is contained in a hyperplane $H_j$, which is in turn a translation of some linear subspace $U_j\subset\R^n$ of dimension $n-1$. If $U_j$ has a basis $\vec{v}_{j,1},\hdots,\vec{v}_{j,n-1}$, then it can be shown that
%\[
%H_j = F_j+\sum_{i=1}^{n-1} \vec{v}_{j,i} \cdot\mathbb{Q}
%\]
%In other words, every element of the hyperplane $H_j$ containing $F_j$ can be written as a sum of an element of $F_j$ and a rational linear combination of the basis vectors of $U_j$.
%Repeating this process for each of the $m$ faces of $\P$, we let
%\[B' = B + \sum_{j\in J} \sum_{i=1}^{n-1} \vec{v}_{j,i}\cdot\mathbb{Q}\]
%be the result of replacing each polytope face in $B$ with the entire hyperplane containing it. As each set $\vec{v}_{j,i}\cdot\mathbb{Q}$ is countable, we have $\dim_P B'=\dim_P B$.
By the countable stability of packing dimension, $\dim_P B' = \dim_P B$.
In particular, $B'$ contains the sets 
\[
A'_{j} = \bigcup_{a \in A_j} U_j + a \cdot \vec{u}_j \subset \R^n.
\]
So,
\[
\dim_P B = \dim_P B' \geq \dim_P \bigcup_{j\in J}A'_j = n-1 + \max_{j\in J} \dim_P A_j.
\]
by finite stability of fractal dimension.
Finally, set
\[
A = \bigcup_{j\in J}\frac{1}{d_j} \cdot A_j = \Big\{ j \cdot x + r(x) : x \in [0,1], j \in J \Big\} \subset \R.
\]
By construction, $A$ contains an $m$-term arithmetic progression with common difference $x$ for every $x \in [0,1]$, and $\dim_P A = \max_{j=1, \dots, m} \dim_P A_j$.
This yields the inequality $\h_n^{P}(\P) \geq n-1 + \f^{P}(m)$.

Conversely, we claim that $\h_n^M(\P) \leq n-1 + \f^{M}(m)$. The argument above can essentially be reversed, with the subtlety that Minkowski dimension is not defined for unbounded sets, and is only finitely stable. Let $A \subset \R$ be a bounded set containing an $m$-term arithmetic progression starting at some $a(x) \in \R$ for every common difference $x \in [0,1]$. Translating by an appropriate real number, we can assume that $a(x) > 0$ for every $x \in [0,1]$.
Set $r(x) = a(x)$, and let $B$ be the corresponding union of homothets of $\P$. Note that $B$ is bounded because $A$ is. As in the other direction of the proof, we proceed to define a larger set $B'$ by modifying $B$.

The intuition is that instead of replacing each face by an entire hyperplane, we replace it by a larger---yet bounded---piece of a hyperplane, taken to be large enough that $B'$ contains $B$. More precisely, with the notation above, $B'$ is composed by a collection of sets $A'_j$, one for each face $j$. But instead of containing translations of $(n-1)$-dimensional linear subspaces $U_j$, $A'$ is taken to contain translations of sufficiently large $(n-1)$-cubes $C_j\subset U_j$.

As before, we have $\dim_M B\leq \dim_M B' = n-1+\dim_M A$. In conclusion, we have the inequalities
\[
n-1 + \f^{P}(m) \leq \h_n^{P}(\P) \leq \h_n^{M}(\P) \leq n-1 + \f^{M}(m),
\]
and by \autoref{thm:AKCequivs}, we have equality everywhere.
\end{proof}
\begin{remark}
In the theorem statement, $\P$ is implicitly an $m$-tope (has $m$ faces). The proof in fact allows for a generalization to polytopes with more faces. Suppose that $\P$ has bounding hyperplanes $\{H_j^{(i)} : j\in J, i\in I\}$ for $J$ an arithmetic progression of length $m$ and $I$ a secondary indexing set intended to allow for additional faces. As usual, each hyperplane $\smash{H_j^{(i)}}$ is determined by the pair $\smash{(\vec{u}_J^{(i)},d_j^{(i)})}$. If $\smash{jd^{(i)}_j = (\vec{u}_j^{(i)})_1}$ for all $i,j$, then the conclusions of \autoref{thm:harmonicpolytopestrong} hold. To see this, note that the sets $\smash{A_j^{(i)}}$ are entirely determined by the first coordinate of $\smash{\vec{u}_j^{(i)}}$ and $\smash{d_j^{(i)}}$, and the scaling factors $r(x)$. Then the subsets of $B'$ containing homothets of $\smash{H^{(i)}_j}$ and $\smash{H^{(i')}_j}$ map identically into $A$, for any $i,i'$. Since $\P$ has finitely many faces, finite stability of Minkowski and packing dimensions is enough to conclude $\h_{n}^{M}(\P) = \h_n^{P}(\P) = n-1 + \f(m)$.
Note that this bound gets worse as we increase the size of $I$, as $m$ gets correspondingly smaller compared to the number of faces of $\P$. 
\end{remark}

\begin{figure}
    \centering
    \includegraphics[width=0.6\textwidth]{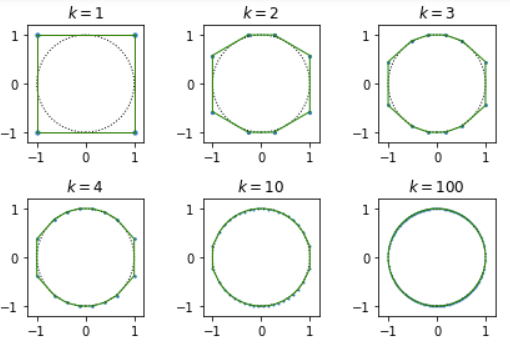}
    \caption{Convergence of the polygons $\P_k$ to the circle when $n=2$.}
    \label{fig:polygons}
\end{figure}

To prove \autoref{thm:harmonicpolytope}, we apply the previous theorem with an appropriate sequence of polytopes $(\P_k)_{k}$.
\begin{proof}[Proof of \autoref{thm:harmonicpolytope}]
For each $k\geq1$, the strategy is to apply the remark following \autoref{thm:harmonicpolytopestrong} with $J_k=\{-k,\dots,k\}$, and $I$ as large as necessary to complete the polytope $\P_k$. Let
\begin{align*}
\vec{u}^{(i)}_{j,k} &= (j/k, \vec{v}^{(i)}_{j,k}),\\
d^{(i)}_{j,k} &= 1,
\end{align*}
for each $j \in J_k=\{-k, \dots, k\}$, where $\vec{v}^{(i)}_{j,k} \in \R^{n-1}$ is chosen so that $||\vec{u}^{(i)}_{j,k}|| =1$, and so that the corresponding polytopes converge to the sphere in the Hausdorff metric.
%and for each $k$, the points $\{\vec{u}^{(i)}_{j,k}\}_{i,j}$ form an $o_k(1)$-net of $S^{n-1}$. By construction, this sequence of polytopes converges to $S^{n-1}$.
To conclude, \autoref{thm:harmonicpolytopestrong} and \autoref{thm:1dAKCequivs} imply that
\[
\lim_{k\to\infty} \h_n^M(\P_k) = \lim_{k\to\infty} \h_n^P(\P_k) = (n-1) + \lim_{k\to\infty} \lim_{N\to\infty} \frac{\log F_k(N)}{\log N}
\]
\end{proof}
For clarity, we illustrate the proof method for $n=2$.\vspace{0.6em}
\begin{example}\label{ex:2D}
When $n=2$, for each fixed $j,k$, there are exactly 2 possible choices of $\smash{v_{j,k}^{(i)}}$: they are $\smash{\pm\sqrt{1-j^2/k^2}}$. Letting $I=\{+,-\}$ parametrize these two possibilities, the collection $\smash{\{H_j^{(+)}:j\in J\}}$ can be thought of a set of bounding hyperplanes for the top half of the polygon $\P_k$, and $\smash{\{H_j^{(-)}:j\in J\}}$ as a set of bounding hyperplanes for the bottom half. The polygons $\P_k$ and their convergence to the circle are illustrated in \autoref{fig:polygons}.
\end{example}
The construction above yields a polygon $\P_1$ with four sides defined by the unit normal vectors $(1,0)$, $(0,-1)$, $(-1,0)$, $(0,1)$, each side a unit distance from the origin. Then as in \autoref{fig:polygons}, $\P_1$ is a square centered at the origin. The above theorem therefore provides information on the minimal Minkowski and packing dimensions of a set containing an axis-aligned squared centered at every point in the straight line segment $[0,1]\times\{0\}$. 
\begin{corollary}\label{cor:squarebounds}
The minimal Minkowski and packing dimensions of a set containing an axis-aligned square centered at each point on a segment of the $x$-axis are equal, and they are between 17/11 and 7/4. That is,
\[
1.\overline{54}\leq \h_2^M(\P_1) = \h_2^P(\P_1) \leq 1.75
\]
\end{corollary}
\begin{proof}
The Katz-Tao inequality asserts that $F'_3(N)\geq N^{6/11}$ \cite{KatzTao}. Applying \autoref{thm:AKCequivs}, we have $\mathcal{f}(3)\geq 6/11$. The lower bound then follows immediately from \autoref{thm:harmonicpolytopestrong} applied to $\P_1$. On the other hand, in Theorem 1.2 of \cite{Squares} it was shown that there is a compact set of Minkowski and packing dimension $7/4$ containing an axis-aligned square with every center in $[0,1]^2$. After translating and scaling if necessary, this set has property specified in the corollary statement.
\end{proof}
Above we consider our set of centers $S$ to be an axis-parallel line segment. Perhaps surprisingly, rotating $S$ by $45^\circ$, we obtain a different answer.
\begin{corollary}\label{cor:diamond}
The minimal Minkowski and packing dimensions of a set containing an axis-aligned square centered at each point on a segment of the line $y=x$ is 1.5.
\end{corollary}
\begin{proof}
Rotating by $45^\circ$, it suffices to consider sets containing a ``diamond'' (a $45^\circ$ rotation of an axis-aligned square) at each point on $[0,1]\times\{0\}$. Such a diamond $\P$ is defined by the four lines determined by the four normal vectors $(1/\sqrt{2},\pm1/\sqrt{2})$ and $(-1/\sqrt{2},\pm1/\sqrt{2})$, and the unit distances $d_j^\pm=1$. Here $J=\{-1\sqrt{2},1/\sqrt{2}\}$ and $I=\{+,-\}$ (as in \autoref{ex:2D}). By the remark following \autoref{thm:harmonicpolytopestrong}, $\h_2^P(\P)=\h_2^M(\P)=1+\f(2)$, so it suffices to compute $\f(2)$. A trivial sumset bound shows that $\f(2)\geq1/2$. On the other hand, it was shown in section 2.3 that $G_2'=O(N^{1/2})$, so it follows from \autoref{thm:AMPimpliesAKC} and \autoref{thm:AKCequivs} that $\f(2)\leq 1/2$, concluding the proof.
\end{proof}

\begin{remark}
The main distinction between the previous two corollaries is that differing slopes of $S$ results in different lengths of the arithmetic progression $J$. In \autoref{cor:diamond}, $|J|=2$; while in \autoref{cor:squarebounds}, $|J|=3$. By appropriately choosing a segment $S$ with slope $1/3$, we achieve $|J|=4$. In this case, one can follow the proof of \autoref{cor:squarebounds} to achieve the (slightly better) lower bound 11/7, using the 4-term Katz-Tao inequality from \cite{KatzTao}.
\end{remark}

So far, we have considered only the case when our set $S$ of centers is a straight line. Not much is known about other choices of $S$. One natural direction of generalization would be to other sets $S$ of dimension 1. In \cite{Squares}, a related construction was given:
\begin{theorem}
There exist sets $S,B\subset\R^2$ such that $B$ contains a square centered at every point of $S$, and
$$\dim_H S=\dim_M S=\dim_P S=1$$
while
$$\dim_M B=\dim_P B = 1+\frac{3}{8}.$$
\end{theorem}

\subsection{Bounds for homothets of polytopes}
In contrast to the previous section, here we consider the case when $B$ contains a homothet of a polytope $\P$ centered at each point of $S=[0,1]^n$.
First, we show a lower bound for such polytope problems in arbitrary dimensions, and then show that this is sharp for simplices. 

\begin{theorem}\label{thm:polytopes}
Let $\P \subset \R^n$ be a convex polytope. 
If $B \subset \R^n$ contains a homothetic image of $\P$ centered at each point in $S \subset \R^n$, then
\[
\dim_P B \geq n-1 + \frac{1}{n+1}\dim_P S.
\]
\end{theorem}
%In particular, when $S = [0,1]^n$, \autoref{thm:polytopes} is sharp for simplices.
\begin{proof}
Our proof is similar to that for the special case of triangles, outlined in the concluding remarks of \cite{Squares}. 
Suppose $B$ contains a homothet of $\P$ $x + r(x) \cdot \P$ for every $x \in [0,1]^n$.
%Suppose that $P$ has $k$ $(n-1)$-dimensional faces, where the $j$\textsuperscript{th} face is specified by the pair $(\hat{u}_j, d_j)$ as outlined in the preliminary section. Further, let $H_j$ denote the $(n-1)$-dimensional hyperplane containing the $j$\textsuperscript{th} face.
Fix an edge $e$ of $\P$, and the corresponding line $\ell$.
Let $v_1$ and $v_2$ be the two vertices of $e$.
Any homothet of $\P$ is uniquely determined by the images of $v_1$ and $v_2$.
We know that $e$ is the intersection of the boundary of $\P$ with $n-1$ hyperplanes---say, $H_{1}, \dots, H_{n-1}$.
There exist hyperplanes $H_{n}$ and $H_{n+1}$ so that
\[
\bigcap_{i=1}^{n-1} H_{i} \cap H_{n} = v_1, \quad \bigcap_{i=1}^{n-1} H_{i} \cap H_{n+1} = v_2.
\]
As in the proof of \autoref{thm:harmonicpolytopestrong}, we parametrize each homothet of $\P$ by the distances of the translates of $H_1, \dots, H_{n+1}$ from the origin.
That is, with $(\vec{u}_i, d_i)$ as earlier, set
\[
A_i = \Big\{ x \cdot \vec{u}_i + (r(x)-1) \cdot d_i: x \in [0,1]^n \Big\}, \quad A'_i = \bigcup_{a \in A_i} H_i + a, \quad i=1, \dots, n+1.
\]
By the same arguments as earlier, we have on one hand that
\[
\dim_P B = n-1 + \max_{i = 1, \dots, n+1} \dim_P A_i,
\]
and on the other
\[
\dim_P S \leq (n+1) \max_{i = 1, \dots, n+1} \dim_P A_i.
\]
Putting this together,
\[
\dim_P B \geq n-1 + \frac{1}{n+1}\dim_P S.
\]

\end{proof}

For a complementary upper bound, we apply probabilistic results of Shmerkin; the following lemma is a modification of corollaries 5.8 and 5.9 of \cite{Pablo}.

\begin{lemma}\label{lem:pablo}
Let $\Gamma$ be a collection of lines in $\R^m$.
For each line $V \in \Gamma$, let $\tilde{V}$ denote the unique line in $\R^m$ parallel to $V$ and passing through the origin.
Suppose that,
\begin{enumerate}[(a)]
    \item for each $V \in \Gamma$ and $I \subsetneq [m]$, $\tilde{V}$ is not contained in the coordinate hyperplane $H^{I} := \{ x \in \R^m: x_i = 0, \forall i \in I \}$, and
    \item there is a cube $Q$ disjoint from the union of diagonals $\Delta = \{ x \in \R^m: x_i = x_j \text{ for some }i \neq j \}$ such that $int(Q) \cap V \neq \emptyset$ for each $V \in \Gamma$.
\end{enumerate}
Then, for every $\epsilon > 0$, there is a suitable set $A \subset \R$ such that $\dim_P A = \dim_M A = 1 - 1/m + \epsilon$ and $A^m \cap V \neq \emptyset$ for each $V \in \Gamma$.
\end{lemma}
%This lemma is a modification of corollaries 5.8 and 5.9 of \cite{Pablo}, and will be explained in further detail in the appendix.

\begin{theorem}\label{thm:simplex}
Let $\P$ be an $n$-dimensional polytope with $m$ faces, i.e.\ $m$ hyperlanes bounding it.
For every $\epsilon > 0$, there is a set $B = B_{\epsilon} \subset \R^n$ containing a homothetic copy of $\mathcal{P}$ centered at every point $x \in \R^n$ such that
\[
\dim_P B = n- \frac{1}{m}+\epsilon.
\]
\end{theorem}
\begin{proof}
Let $(\vec{u}_i, d_i)$, $1 \leq i \leq m$ parametrize the $m$ faces of $\P$.
For each $x \in \R^n$, define
\[
x \cdot U = (\langle x, \vec{u}_1 \rangle, \dots, \langle x, \vec{u}_m \rangle) \in \R^m,
\]
and
\[
V_x = \{ x \cdot U + r (d_1, \dots, d_m): r \in \R \} \subset \R^m.
\]
Note that each $d_i$ is nonzero, so for each $V_x$, the corresponding $\tilde{V}_x$ is not contained in any of the coordinate hyperplanes.
We want to find a cube $S \subset \R^n$ so we can apply the lemma to $\Gamma = \{ V_x \}_{x \in S}$.
$V_x$ intersects $\Delta$ if and only if $\langle x , u_i \rangle + r \cdot d_i = \langle x, u_j \rangle + r \cdot d_j$ for some pair $i \neq j$.
If $d_i = d_j$, then this equation is satisfied if and only if $x$ lies on the hyperplane orthogonal to $u_i - u_j$.
Since there are only finitely many such ``bad'' hyperplanes, we can choose a cube $S$ disjoint from them.
If $d_i \neq d_j$, since this equation is a continuous function of $x \in S$, there is some threshold $\delta$ such that for all $r > \delta$, this equation is not satisfied.
In other words, the sets
\[
V_x' := \{ x \cdot U + r (d_1, \dots, d_m): r > \delta \}, x \in S
\]
are disjoint from $\Delta$.
If a cube $Q$ disjoint from $\Delta$ intersects some $V_w'$ in its interior, then there is some neighbourhood $O$ of $w$ such that $int(Q)$ intersects $V_x'$ for all $x \in O$.
Shrinking $S$ is necessary, we assume that $S \subset O$ so that (b) holds.
Applying the lemma, there is a set $A \subset \R$ such that for each $x \in \R^n$, $A^m \cap V_x \neq \emptyset$, or, $A$ contains points $\langle x,\vec{u}_1 \rangle + r\cdot d_1, \dots, \langle x,\vec{u}_m \rangle + r\cdot d_m$ for some scaling factor $r$.

Now, defining
$A_i := \bigcup_{a \in A}H_i + (a-d_i) \cdot \vec{u_i}$, define $B = \bigcup_{i=1}^{m}A_i$, so that $B$ contains a homothet of $\P$ centered at each point of $S$, and $\dim_P B = n-  \frac{1}{m} + \epsilon$.
\end{proof}

\section{Patterns in other fields}\label{sec:fields}
%\subsection{A detour}
%We use \autoref{lem:pablo} to obtain an upper bound for $G_{U_k}'(N)$, for any set $U_k \subset \Z\setminus\{0\}$ of size $k$.
%This bound will be improved later when we restrict to sets $U_k$ that grow subexponentially in $k$.
%\begin{theorem}
%$\g^P(U_k) = \g^{M}(U_k)\leq 1-1/k$.
%\end{theorem}
%\begin{proof}
%Let $U_k = \{u_1, \dots, u_k\}$.
%For $x \in [0,1]$, define
%\[
%V_x = \big\{ (x, \dots, x) + r \cdot (u_1, \dots, u_k): r \in \R \big\} \subset \R^k.
%\]
%Since each $u_i$ is nonzero, for each $V_x$, the corresponding $\tilde{V}_x$ is not contained in any of the hyperplanes.
%Further, since $u_1, \dots, u_k$ are all distinct, $V_x$ only intersects $\Delta$ in the point $(x, \dots, x)$, so we can find cubes $Q_x$ satisfying (b).
%Applying \autoref{lem:pablo} gives a set $A \subset \R$ with $\dim_P A = \dim_M A = 1 - 1/k + \epsilon$ such that for each $x \in [0,1]$, $A^k \cap V_x \neq \emptyset$.
%Equivalently, $A$ contains the points $x + r \cdot u_1, \dots, x + r \cdot u_k$ for some $r \in \R$, so $A$ contains a $U_k$-pattern for every basepoint in $[0,1]$.
%\end{proof}
%
%\begin{corollary}
%For any $\epsilon > 0$, there is $N$ sufficiently large so that $G'_{U_k}(N) \leq N^{1 - 1/k + \epsilon}$.
%\end{corollary}

\subsection{Finite fields}
Given a prime $p$, let $\pi_p: \Z \to \F_p$ be the canonical projection.
\begin{definition}
Given $U \subset \Z$, a \textit{$U$-pattern} in $\F_p^n$ is a set of the form $x + \pi_p(U) \cdot r$ for $x,r \in \F_p^n$, with $r$ nonzero.
\end{definition}
The following conjecture is technically a family of conjectures for each $n \in \N$ and each sequence $(U_k)_{k \in \N}$ such that $|U_k| = k$.
\begin{conjecture}[The $n$-dimensional finite field $(U_k)$-pattern conjecture]
Let $g_{n,U_k}(p)$ be the size of the smallest set that contains a $U_k$-pattern with basepoint $x$, for every $x \in \F_p^n$.
Then,
\[
\lim_{k \to \infty} \lim_{p \to \infty} \frac{\log g_{n,U_k}(p)}{\log p} = n.
\]
\end{conjecture}
Our first main result is the following.
\begin{theorem}
For any sequence $(U_k)_{k \in \N}$, the pattern conjecture in $\Z$ is equivalent to the  $n$-dimensional finite field $(U_k)$-pattern conjecture for every $n \in \N$.
\end{theorem}
The proof proceeds in three steps.
\begin{step}[1]\label{prop:ZimpliesFpn}
The $(U_k)$-pattern conjecture in $\Z$ implies the  $n$-dimensional finite field $(U_k)$-pattern conjecture for every $n$.
\end{step}
\begin{step}[2]\label{prop:FpnimpliesFp}
For any $n \in \N$, the  $n$-dimensional finite field $(U_k)$-pattern conjecture implies the $1$-dimensional version.
\end{step}
\begin{step}[3]\label{prop:FpimpliesZ}
The  $1$-dimensional finite field $(U_k)$-pattern conjecture implies the $(U_k)$-pattern conjecture in $\Z$.
\end{step}

For the third step, we need the following lemma, which follows from \autoref{lem:pablo}.
\begin{lemma}\label{lem:GUkbound}
For any set $U_k \subset \N$ of size $k$, $\g^P(U_k) = \g^{M}(U_k)\leq 1-1/k$.
As a result, for any $\epsilon > 0$, there is $N$ sufficiently large so that $G'_{U_k}(N) \leq N^{1 - 1/k + \epsilon}$.
\end{lemma}
\begin{proof}
Let $U_k = \{u_1, \dots, u_k\}$.
For $x \in [0,1]$, define
\[
V_x = \big\{ (x, \dots, x) + r \cdot (u_1, \dots, u_k): r \in \R \big\} \subset \R^k.
\]
Since each $u_i$ is nonzero, for each $V_x$, the corresponding $\tilde{V}_x$ is not contained in any of the hyperplanes.
Further, since $u_1, \dots, u_k$ are all distinct, $V_x$ only intersects $\Delta$ in the point $(x, \dots, x)$, so we can find a cube $Q$ satisfying (b).
Applying \autoref{lem:pablo} gives a set $A \subset \R$ with $\dim_P A = \dim_M A = 1 - 1/k + \epsilon$ such that for each $x \in [0,1]$, $A^k \cap V_x \neq \emptyset$.
Equivalently, $A$ contains the points $x + r \cdot u_1, \dots, x + r \cdot u_k$ for some $r \in \R$, so $A$ contains a $U_k$-pattern for every basepoint in $[0,1]$.
\end{proof}

\begin{proof}[Proof of \href{prop:ZimpliesFpn}{step 1}]
Fix a set $U_k$ and $n \in \N$.
Fix $\epsilon > 0$, and let $p$ be a large enough prime so that $|u| < p^{\epsilon}$ for every $u \in U_k$.
Suppose $A_1 \subset \F_p^n$ contains a $U_k$-pattern for every basepoint $x \in \F_p^n$.
Let $\psi: \F_p \to \{0, \dots, p-1 \}$ be the natural ``lifting'' map, and extend this to a map $\psi_n: \F_p^n \to \{0, \dots, p-1\}^n$ by ``lifting'' each coordinate.
For each $x \in \F_p^n$, we have a pattern $x + r(x) \cdot U_k$ in $A_1$.
Define
\[
A_2 = \Big\{ \psi_n(x) + \psi_n(r(x)) \cdot u : x \in \F_p^n, u \in U_k \Big\} .
\]
By construction, $A_2 \subset [-p^{\epsilon}(p-1), p^{\epsilon}(p-1)]^{n}$, and $\pi(A_2) \subset A_1$, where $\pi: \Z^n \to \F_p^n$ is the natural projection.
Each coordinate of an element in $\pi(A_2)$ has at most $2p^{\epsilon}$ possible preimages, so $|A_2| \leq \big( 2p^{\epsilon} \big)^n |A_1|$.
Next, choose $b \in \N$ large enough so that the following map is injective.
Further, if $b$ is sufficiently large, then
\[
\phi \Big( (x_1, \dots, x_n) + u \cdot (r_1, \dots, r_n) \Big) = \phi \Big( (x_1, \dots, x_n) \Big) + u \cdot \phi \Big( (r_1, \dots, r_n) \Big).
\]
In other words, $\phi(A_2) \subset \Z$ contains a $U_k$-pattern for $p^n$ different basepoints, so $|A_2| \geq G'_{U_k}(p^n)$.
This yields $g_{n,U_k}(p) \geq (2p^{\epsilon})^{-n} G_{U_k}(p^n)$, so
\[
\lim_{p \to \infty} \frac{\log g_{n,U_k}(p)}{\log p} \geq n \lim_{p \to \infty} \frac{\log G'_{U_k}(p^n)}{\log p^n} - n \epsilon.
\]
Since $n$ is fixed and $\epsilon > 0$ is arbitrary,
\[
\lim_{p \to \infty} \frac{\log g_{n,U_k}(p)}{\log p} \geq n \lim_{p \to \infty} \frac{\log G'_{U_k}(p^n)}{\log p^n}.
\]
\end{proof}

\begin{proof}[Proof of \href{prop:FpnimpliesFp}{step 2}]
The proof is straightforward. Suppose $A \subset \F_p$ contains a $U_k$ pattern for each basepoint $x \in \F_p$.
Then, $A^n \subset \F_p^n$ contains a $U_k$-pattern for each basepoint $(x_1, \dots, x_n) \in \F_p^n$.
This yields $g_{n,U_k}(p) \leq \Big( g_{1,U_k}(p) \Big)^n$, so
\[
\lim_{p \to \infty} \frac{\log g_{1,U_k}(p)}{\log p} \geq \frac{1}{n} \lim_{p \to \infty} \frac{\log g_{n,U_k}(p)}{\log p}.
\]
\end{proof}

\begin{proof}[Proof of \href{prop:FpimpliesZ}{step 3}]
Choose $0 < \epsilon 1/k$, and let $N$ be as in \autoref{lem:GUkbound}.
Let $A_1$ be a set of size $G_{U_k}(N)$ containing a $U_k$-pattern for every basepoint in $[N]$, so that
\[
|A_1| \leq C_k \log N \cdot N^{1 - \epsilon}.
\]
We want a prime $p > N$ so that projecting $A_1$ to $\F_p$ gives a $U_k$-pattern with $N$ basepoints in $\F_p$.
However, given a pattern $x + r(x) \cdot U_k$, we need to guarantee that $r(x)$ is nonzero modulo $p$.
Simultaneously, we need to guarantee that $p$ is close to $N$, so that $g_{1,k}(p)$ is close to $G_{U_k}(p)$.
For this second condition, we want to choose $p$ to be a prime in the interval $(N, 2N]$.

Suppose we have no such choice: for every prime $p \in (N, 2N]$, there exists $x \in [N]$ so that $p$ divides $r(x)$.
This is where we use the minimality of $A_1$.
For fixed $u \in U_k$, we must have
\[
\Big| \{ x + r(x) \cdot u: x \in [N] \} \Big| \leq C_k \log N \cdot N^{1 - \epsilon}.
\]
In particular, we must have ``many'' nontrivial solutions to the equation
\[
x + r(x) \cdot u = y + r(y) \cdot u: \quad x,y \in [N].
\]
We will show that there exists a large subset $S \subset [N]$, such that this equation does \textit{not} hold for any pair $x,y \in S$.
The multiset $\{r(x) : x \in [N] \}$ has cardinality $N$, while there are only $N/\log N$ primes $p \in (N, 2N]$.
As a result, there exists a prime $p \in (N, 2N]$ such that for at least $N - N/\log N$ values of $x \in [N]$,  $p$ divides $r(x)$.
Let $S \subset [N]$ be the set of these elements $x \in [N]$.
Suppose for two distinct elements $x,y \in S$, $u \cdot \big( r(y) - r(x) \big) = (x-y)$.
Then $r(x)$ and $r(y)$ must be distinct, and $p$ divides $r(y) - r(x)$, but $p$ cannot divide $x-y$ because $p>N$.
So, no such solution exists.
This yields that $|A_1| \geq N - N/\log N$, but this contradicts the upper bound for large enough $N$.

In other words, we can choose a prime $p=p(N) \in (N, 2N]$ such that $r(x) \neq 0 \Mod{p}$ for any $x \in [N]$.
Then, $A_2 = \pi_p(A_1)$ contains a $U_k$-pattern for $N \geq p/2$ distinct basepoints in $\F_p$.
By translating $A_2$ by some set $T$, we obtain a set $A_3$ such that $A_3$ contains a $U_k$-pattern for every basepoint in $\F_p$, and $|A_3| \ll \log p \cdot |A_2| \leq \log p \cdot |A_1|$.
So,
\[
G_{U_k}(N) \geq \frac{1}{\log p}g_{1,U_k}(p(N)).
\]
It follows that
\[
\lim_{N \to \infty} \frac{\log G_{U_k}(N)}{\log N} \geq \lim_{N \to \infty} \frac{\log g_{1, U_k}(p(N))}{\log p(N)}.
\]
\end{proof}

\subsection{An upper bound for number fields}

In \cite{GreenRuzsa}, Green and Ruzsa show the upper bound $F_{k}(N) \leq N^{1 - c/\log \log k}$.
We use a similar technique to obtain an upper bound for patterns in $\Z^d$.

For the rest of this section, we fix the following notation.
Let $n$ be an integer of the form $2,4,p^k$ or $2p^k$ for an odd prime $p$, and $d = \phi(n)$, where $\phi$ denotes the Euler totient function.
Let $\zeta$ be a primitive $n$th root of unity.
Define the free $\Z$-module
\[
\OK = \Z[\zeta] \cong \Z \oplus \Z \cdot \zeta \oplus \dots \oplus \Z \cdot \zeta^{d-1}.
\]
The latter isomorphism is a well-known fact in number theory (see for example Proposition 10.2 in \cite{Neukirch}).
Define a $\Z$-module isomorphism $\psi: \Z^d \to \OK$ in the natural way.
This allows us to define a multiplication $\otimes$ on $\Z^d$ by $r \otimes u = \psi^{-1}(\psi(r) \cdot \psi(u))$.

\begin{theorem}[Proposition 10.3, \cite{Neukirch}]
Let $a$ be a primitive root modulo $n$, i.e.\ $a$ generates $\Z_{n}^{\times}$.
If $q$ is a prime in $\Z$ and $q \equiv a \Mod{n}$, then $q\OK$ is a prime ideal in $\OK$.
\end{theorem}
\begin{definition}
If $U \subset \N^d$, a $U$-pattern is a set of the form $x + r \otimes U = \{x + r \otimes u: u \in U \}$, for $x, r \in \Z^d$.
We call $x$ and $r$ the \textit{basepoint} and \textit{scaling factor} of the pattern respectively.
\end{definition}
Note that this is different from the kinds of patterns we considered in earlier sections, as the scaling factor also belongs to $\Z^d$.
Let $(U_k)_{k \in \N}$ be a sequence such that each $U_k \in \binom{\mathbb{N}^d}{k}$.\footnote{Recall that by our convention, $0 \notin \N$, so for each $u \in U_k$, no coordinate of $u$ is equal to zero.}
We will assume for convenience that from here on, the $\ell^\infty$ norms of the elements of $U_k$ are bounded by some subexponential $f(k)$, i.e.
\[
\max_{u\in U_k} ||u||_\infty \leq f(k)=2^{o(k)}.
\]

\begin{theorem}\label{thm:numberfields}
Let $U_k\in\binom{\mathbb{Z}^{d}}{k}$ be a sequence as above.
Then for all large $k$ and $N$, there is a compact set $B_k$ containing a $U_k$-arithmetic pattern with every center in $[N]^{d}$ such that
$$|B_k| < N^{d-1/O(\log k)}$$
\end{theorem}
The important special cases are when $d=1$ and $d=2$.
We will need three technical lemmas. 
\begin{lemma}\label{lem:boundedproduct}
Let $r,u\in\Z^d$. Then
\[
||r\otimes u||_\infty \leq 2d \cdot \big(||r||_\infty \cdot ||u||_\infty\big).
\]
\end{lemma}
\begin{proof}
Given $r = (r_1, \dots, r_d)$ and $u = (u_1, \dots, u_d)$, their product is
\[
\psi(r \otimes u) = \psi(r)\cdot\psi(u) = \sum_{0\leq i,j < d} r_i u_j \zeta^{i+j} = \sum_{k=0}^{2(d-1)} a_k \zeta^k.
\]
for some integers $a_k \in \Z$.
We want to express this in the basis $(1, \zeta, \dots, \zeta^{d-1})$ so that we can pull the product back to $\Z^d$ and bound its norm.
Writing
\[
\sum_{k=0}^{2(d-1)}a_k \zeta^k = \sum_{l=0}^{d-1}b_l \zeta^l,
\]
it is clear that 
\[
\sup_{l=0, 
\dots, 2(d-1)} \big|b_l\big| \leq 2 \sup_{k=0, \dots, d-1}\big|a_k\big| \leq 2d \sup_{i=0, \dots, d-1}\big|r_i\big| \cdot \sup_{j=0, \dots, d-1} \big|u_j\big|.
\]
\end{proof}
\begin{lemma}[Chinese remainder theorem, Theorem 3.6 \cite{Neukirch}]\label{lem:CRT}
Let $q_1, \dots, q_k \in \Z$ be primes such that $q_i\OK$ is a prime ideal, and $Q = \prod_{i=1}^{k}q_i$.
Then we have a ring isomorphism
\[
\faktor{\OK}{Q\OK}\cong \bigoplus_i \mathbb{F}_{q_i^d}.
\]
\end{lemma}
\begin{lemma}\label{lem:primorial}
Suppose $\gcd(a,d)=1$, and let $p_1<p_2<\hdots$ denote the sequence of all primes, $q_1<q_2<\hdots$ denote those which are congruent to $a$ modulo $n$, where $a$ is a fixed primitive root modulo $n$ as earlier. 
Then
$$\prod_{i=1}^m q_i = \left(\prod_{i=1}^m p_i\right)^{(1+o(1))} = m^{(1+o(1))m}$$
\end{lemma}
\begin{proof}
According to Poussin, the number of primes $q_i\equiv a\Mod{n}$ less than a number $x>0$ is given by
$$\pi_{n,a}(x) = (1+o(1))\frac{1}{n}\pi(x).$$
Together with the prime number theorem, which states that
$$\pi(x) = (1+o(1)) \frac{x}{\log x}$$
we conclude that the $i$\textsuperscript{th} prime $p_i$ is approximately
$$p_i = (1+o(1)) i\log i$$
and that the $i$\textsuperscript{th} prime $q_i$ congruent to $a$ modulo $n$ is approximately
$$q_i = (1+o(1))ni\log i.$$
We then have
\begin{align*}
    \prod_{i=1}^m q_i &= \prod_{i=1}^m (1+o(1)) ni\log i \\
    &= (1+o(1))^m n^m \prod_{i=1}^m i\log i \\
    &= (1+o(1))^m n^m \prod_{i=1}^m (1+o(1))p_i\\
    &= (1+o(1))^{2m} n^m \prod_{i=1}^m p_i.
\end{align*}
Applying a primorial estimate yields
\begin{align*}
    \prod_{i=1}^m q_i &= (1+o(1))^{2m} n^m \prod_{i=1}^m p_i \\
    &= (1+o(1))^{2m} n^m e^{(1+o(1))m\log m}\\ \\
    &= e^{2m\log(1+o(1))}e^{m\log n} e^{(1+o(1))m\log m}\\
    &= e^{(1+o(1))m\log m}.
\end{align*}
The result follows.
\end{proof}

The proof of \autoref{thm:numberfields} is a routine deduction from the following discrete proposition, which we will prove first.
\begin{proposition}
Let $U_k\in\binom{\N^{d}}{k}$ be a sequence whose $\ell^\infty$ norm is subexponential in $k$. Then for all large $k$ and $N$, there is a compact set $B_k$ containing a $U_k$- pattern with every basepoint in $[N]^{d}$ such that
$$|B_k| < N^{d-1/O(\log k)}$$
\end{proposition}
\begin{proof}
Let
\[
\max_{u\in U_k} ||u||_\infty \leq f(k) = 2^{o(k)}
\]
as earlier.
Let the residue classes of $\mathbb{Z}_n$ which generate $\mathbb{Z}_n^\times$ be denoted by $a_1,\hdots, a_{\phi(d)}$. Denote the sequence of primes which are congruent to some $a$ coprime to $n$ by $q_1<q_2<\hdots$. Fix $m=2f(k)$, and let $Q(k) = Q=\prod_{i=f(k)}^m q_i$.
By \autoref{lem:primorial},
\[
\log Q = (1+o(1))m\log m - (1+o(1)) f(k)\log f(k) = (1+o(1)) f(k)\log f(k).
\]
Let $U_k = \{u_1, \dots, u_k \}$.
For the rest of the proof, we will denote $[A]^d=\{(x_1,\hdots,x_d) \in \mathbb{Z}^d : x_i\in[A]\}$. Additionally, let
$$\psi : \mathbb{Z}^d\longrightarrow\mathcal{O}_K$$
be the aforementioned isomorphism of $\mathbb{Z}$-modules. For convenience, we will denote $\psi(x)=\overline{x}$. Let
$$\pi_i : \mathcal{O}_K\longrightarrow \faktor{\mathcal{O}_K}{(q_i)}\cong \mathbb{F}_{q_i^d}$$
be the natural projection maps. The map $\pi = \prod_i \pi_i$ is then the projection  onto $\mathcal{O}_K/(Q)$.

Now, for each $x\in[Q-1]^d$, consider the projections $\pi_i(\overline{x})$ onto each of the finite fields $\mathbb{F}_{q_i^d}$. 
By the Chinese remainder theorem, there is a unique element $r(x)\in[Q]^d$ such that
\[
\pi_i(\overline{r(x)}) = \pi_i(\overline{x})^2
\]
for each $i= f(k), \dots, m$. 
Define the set
\[
S = \{x+u_j \otimes r(x) : j \in [k], x\in[Q-1]^d\}.
\]
That is, $S$ contains a $U_k$-pattern with basepoint $x$ and scaling factor $r(x)$ for each $x \in [Q-1]^d$.
We want to show that $S$ is small by showing that each projection to $\F_{q_i^d}$ is small, and that the preimage of each point of $S$ under this projection is also small.

Fix some $j\in[k]$. 
For each $i$,
\[
\pi_i\Big(\overline{x+u_j\otimes r(x)}\Big) =  \pi_i(\overline{x}) + \pi_i(\overline{u_j})\cdot\pi_i(\overline{x})^2 = \pi_i(\overline{u_j})\left(\pi_i(\overline{x})+\frac{1}{2\pi_i(\overline{u_j})}\right)^2 - \frac{1}{4\pi_i(\overline{u_j})}
\]
in $\mathbb{F}_{q_i^d}$. Note that since each coordinate of $u_j$ is nonzero and at most $f(k) < q_i$, $\pi_i(\overline{u_j})$ is invertible in $\F_{q_i^d}$, so the quantity on the right is well-defined.
As $x$ varies, these quantities are in bijection with the quadratic residues in $\mathbb{F}_{q_i^d}$, and so $\pi_i(\overline{x+u_j\otimes r(x)})$ takes values in a set of size at most $(q_i^d+1)/2$ in $\mathbb{F}_{q_i^d}$.
By the Chinese remainder theorem, $\pi(\overline{x + u_j\otimes r(x)})\in\mathcal{O}_K/(Q)\cong\mathbb{Z}^d/Q\mathbb{Z}^d$ takes values in a set of size at most $2^{-f(k)} \prod_{i=f(k)}^m (q_i^d+1)$. 
Lifting to $\Z^d$, and applying \autoref{lem:boundedproduct}, we have 
\[
||x+u_j\otimes r(x)||_\infty \leq ||x||_\infty + O\big(||r(x)||_\infty ||u_j||_\infty\big) \leq Q + O(Qf(k)) = O(Qf(k))
\]
since $x,r(x)\in[Q]^d$. 
Therefore, this quantity takes values in a set of size $O(f(k)^{d})\cdot 2^{-f(k)} \prod_{i=f(k)}^m (q_i^d+1)$.
Summing over all values of $j$,
\[
|S| \leq kO(f(k)^{d}) 2^{-f(k)}Q^{d}\prod_{i=f(k)}^m \Big(1+\frac{1}{q_i^d}\Big)
\]
First, we apply \autoref{lem:primorial} to obtain the bound $\prod_{i=f(k)}^m (1+1/q_i^d) < \prod_{i=1}^m (1+1/q_i) \ll\log m \ll k$.
Next, an elementary counting argument\footnote{Each element of $U_k$ has $\ell^{\infty}$ norm at most $f(k)$.} tells us that $f(k)^d \geq k$, so
\begin{align*}
    |S| &\ll k^2 f(k)^{d} 2^{-f(k)} Q^{d} \\
    &= 2^{-f(k)(1+o(1))} Q^{d} \\
    &= Q^{d-\frac{1}{O(\log k)}}.
\end{align*}

Since we have only constructed the set $S$ for one integer $Q$, the rest of the proof is a technicality to extend this construction for all integers.
First, $q$ an arbitrary positive integer, let $N_q:= Q^q$ and consider the set
\[
A_q := \{ s_0 + s_1 Q + \hdots + s_{q-1}Q^{q-1} : s_i\in S\}.
\]
Then $|A_q| \leq |S|^q\leq N_q^{d-\frac{c}{\log k}}$. The set $A_q$ contains a $U_k$-pattern for every basepoint $x\in \{0,1,\hdots,Q-1\}^d$, or in other words for all $x\in\{0,\hdots,N_q-1\}^d$.
Finally suppose $N$ is an arbitrarily large positive integer. With $q$ minimal such that $N_q > N$, set $A:=A_q$. Then $A$ contains a $U_k$-pattern for every basepoint $x\in[N]^d$. Moreover,
\[
|A| \ll_k N^{d-\frac{c}{\log k}}
\]
and the result follows.
\end{proof}
We now prove \autoref{thm:numberfields}

\begin{proof}[Proof of \autoref{thm:numberfields}]
Let $Q \in \Z$ and $S \subset \Z^d$ be as earlier.
Define the set 
\[
F = \Big\{\sum_{i=0}^{\infty}s_iq^{-i}: s_i \in S \Big\}.
\]
It follows easily from the definition of Minkowski dimension that  $\dim_M F = d - c/\log k$.
Further, the same argument used in the proof of \autoref{thm:AKCequivs} tells us that $F$ contains a $U_k$-pattern for every basepoint $x \in [0,1]^d$.
\end{proof}

When $n=2$, then $d=1$, and the isomorphism $\Z \to \OK$ is just the identity $\Z \to \Z$.
In this case, we obtain the corollary
\begin{corollary}
If $(U_k)_{k \in \N}$ is a sequence in $\Z$ that grows subexponentially in $k$, then 
\[
\lim_{N \to \infty}\frac{\log G'_{U_k}(N)}{\log N} \leq 1 - c/\log k.
\]
\end{corollary}
Setting $U_k = [k]$ proves \autoref{prop:Gkbound}.

When $n=4$, then $d=2$, and the isomorphism $\Z^2 \to \OK$ is the natural map $\Z^2 \to \Z[i]$, so that the multiplication $\otimes$ is just complex multiplication.
Geometrically, this corresponds to scaling and rotating the sets $U_k$.
Setting $U_k$ to be the vertex set of a lattice $k$-gon, we obtain the proof of \autoref{prop:latticekgon}.
%\todo{Do we need a conclusion?}
\appendix
\section{Appendix}
\subsection{Random translates}
We recall some standard results on random translates as presented in \cite{GreenRuzsa}.
\begin{lemma}
If $S \subset \{1, \dots, X\}$, then there is a set $T$ of size $\ll \frac{X}{|S|} \log X$ such that $S+T \supset \{1, \dots X\}$.
\end{lemma}

\begin{lemma}
If $S \subset \F_p^n$, then there is a set $T \subset \F_p^n$ of size $\ll \frac{p^n}{|S|}n\log p$ such that $S+T = \F_p^n$.
\end{lemma}

%\todo{Should we repeat the proof?}
\subsection{Equivalence of the pattern and restricted pattern conjectures}
Now, we will prove \autoref{prop:restrictedequiv} by modifying the proof from \cite{GreenRuzsa}.
\begin{proof}[Proof of \autoref{prop:restrictedequiv}]

It is clear that $G'_{U_k}(N) \leq G_{U_k}(N)$.
For the converse, we claim that
\[
G_{U_k}(N) \leq C_k \log N \cdot G'_{U_k}(N)
\]
where $C_k$ is a constant depending only on $U_k$.
It will then follow that
\[
\lim_{k \to \infty} \lim_{N \to \infty}\frac{\log G_{U_k}(N)}{\log N} \leq \lim_{k \to \infty} \lim_{N \to \infty}\frac{\log \Big(C_k \log N \cdot G'_{U_k}(N)\Big)}{\log N} = \lim_{k \to \infty} \lim_{N \to \infty}\frac{\log G'_{U_k}(N)}{\log N},
\]
so that the restricted pattern conjecture implies the pattern conjecture.
Let $U_k=\{u_1,\hdots,u_k\}$, with $u_1 < \dots < u_k$. 
Let $A_0$ be a set that contains a $U_k$-term arithmetic progression with distinct integer basepoints $a_1, \dots, a_N$.
Let the corresponding scaling factors be $r_1, \dots, r_N$, not necessarily distinct.
We want to construct a set $A_1$ containing a $U_k$-term progression with basepoints $1, \dots, N$ such that $|A_1| \leq C_k \log N |A_0|$.

Choose $\theta \in (0,1)$ uniformly at random, and define
\[
\phi_{\theta}: \Z \to \{0, \dots, N-1 \}
\]
by
\[
\phi_{\theta}(x) = \lfloor N \{\theta x \}\rfloor.
\]
We want to count the number of pairs $(i,j)$ for which $\phi_{\theta}(a_i) = \phi_{\theta}(a_j)$.
If $i \neq j$, then
\[
\Pr (\phi_{\theta}(a_i) = \phi_{\theta}(a_j)) = \Pr (\lfloor N \{\theta a_i\} \rfloor - \lfloor N \{\theta a_j\} \rfloor = 0)
\leq \Pr (\{\theta |a_i - a_j| \} \leq 1/N) = 2/N,
\]
since $\theta$ was chosen uniformly at random.
The expected number of such pairs is then at most $N-1$, so there is some choice of $\theta$ for which the number of such pairs is at most $N-1$.
For $n \in \{0, \dots, N-1\}$, set
\[
f(n) = \Big| \{i: \phi_{\theta}(a_i) = n \} \Big|.
\]
Then,
\[
\sum_{n} \binom{f(n)}{2} \leq N-1, \quad \sum_{n} f(n) = N.
\]
So, $\sum_{n} f(n)^2 \leq 3N$.
Applying Cauchy-Schwarz,
\[
N^2 = \Big( \sum_{n} f(n) \Big)^2 \leq \Big( \sum_{f(n) \neq 0} 1 \Big)\Big( \sum_{n}f(n)^2 \Big) \leq \Big( \sum_{f(n) \neq 0} 1 \Big) 3N.
\]
It follows that $f(n) \neq 0$ for at least $N/3$ values of $n$, or there are at least $N/3$ distinct values of $\phi_{\theta}(a_i)$.
Now let $A_2 = \phi_{\theta}(A_0)$; $|A_2| \leq |A_0|$, but it is not clear that $A_2$ contains the long progressions we need.
Nevertheless, for every $x,y \in \N$,
\[
\phi_{\theta}(x+y) - \phi_{\theta}(x) - \phi_{\theta}(y) = \lfloor N \{ \theta(x+y) \} \rfloor - \lfloor N \{ \theta x \} \rfloor -  \lfloor N \{ \theta y \} \rfloor \in \{ 0,1 \} - \{0,N \} = \{0,1,-N,1-N \}.
\]
By induction, it follows that
\[
\phi_{\theta}(a_i) + \phi_{\theta}(r_i)\cdot u_j - \phi_{\theta}(a_i + r_i\cdot u_j) \in \{0, \dots, u_j \} - \{0, N, \dots u_j \cdot N \}
\]
for $0 \leq j \leq k-1$.
Define
\[
A_3 = A_2 + \{0, 1, \dots, u_k \} - \{0, N, \dots, u_k \cdot N \}.
\]
Then,
\[
|A_3| \leq (u_k+1)^2 |A_2| \leq (u_k+1)^2 |A_0|.
\]
$A_3$ contains a $k$-term progression centered at at least $N/3$ distinct elements of $\{ 1, \dots, N \}$.
Now, results on random translates tell us that there is a set $T$ of size $C \cdot \log N$ such that every element of $\{1, \dots, N \}$ can be written as $\phi(a_i) + t$.
Then
\[
A_1 = A_3 + T
\]
is the desired set, and
\[
|A_1| \leq C \log N |A_3| \leq C (u_k+1)^2 \log N |A_0|.
\]
Setting $C_k = C (u_k+1)^2$ yields the desired result.
\end{proof}
\subsection{The amplification lemma}
\begin{lemma}[Amplification lemma]
Suppose that there exist $B\subset \Z \times \Z$, a set $S \subset \N$ of size $k$, and $\epsilon>0$ such that
\[
\lvert \pi_{-1}(B) \rvert > \max_{j \in S} \lvert \pi_j(B) \rvert^{1 + \epsilon}.
\]
Then, for all $M>0$, there exists $B'\subset\Z \times \Z$ such that
\[
\lvert \pi_{-1}(B') \rvert > M \max_{j\in S} \lvert \pi_j(B') \rvert^{1 + \epsilon}.
\]
\end{lemma}
\begin{proof}
For each $n>0$, define
\[
B^{\oplus n} = 
\Big\{
\big( (x_1,\hdots x_n),(y_1,\hdots y_n) \big) \in \mathbb{Z}^n\times\mathbb{Z}^n : (x_i,y_i)\in B 
\Big\}
\]
and for each $j \in S$.
\begin{align*}
    \pi_j^{\oplus n} : \mathbb{Z}^n\times\mathbb{Z}^n&\longrightarrow\mathbb{Z}^n\\
    (x,y)&\longmapsto x+jy
\end{align*}
for all $j$. Then we have
\[
\lvert \pi_j^{\oplus n}(B^{\oplus n}) \rvert = \lvert \pi_j(B) \rvert^n
\]
for all $j\in S$, and $n \in \N$. 
Since by assumption $\lvert \pi_{-1}(B) \rvert > \max_{j \in S} \lvert \pi_j(B) \rvert^{1+\epsilon} $, for any $M > 0$, for $n$ sufficiently large,
\[
\lvert \pi_{-1}^{\oplus n}(B^{\oplus n}) \rvert > M \max_{j\in S} \lvert \pi_j^{\oplus n} (B^{\oplus n}) \rvert^{1+\epsilon}.
\]
We will create $B'$ from this set $B^{\oplus n}$. For each $t\in\mathbb{Z}$, consider the maps
\begin{align*}
    \psi_t : \mathbb{Z}^n\times\mathbb{Z}^n&\longrightarrow\mathbb{Z}\times\mathbb{Z} \\
    (x,y)&\longmapsto ((t,\hdots,t^n)\cdot x, (t,\hdots,t^n)\cdot y),
\end{align*}
where the dot denotes the standard inner product.
We define an analogous map
\[
\phi_t: \Z^n \to \Z: \quad \phi_t(x_1, \dots, x_n) = x_1t + x_2t^2 + \dots + x_nt^n.
\]
For fixed $j \in S \cup \{-1\}$, and $(x,y), (x',y')$ in $B^{\oplus n}$, 
\[
\big(\pi_j^{\oplus n}(x,y) - \pi_j^{\oplus n}(x',y')\big)\cdot(t,\hdots,t^n)
\]
is a polynomial in $t$ with at most finitely many roots if $\pi_j^{\oplus n}(x,y) \neq \pi_j^{\oplus n}(x',y')$.
Since $B$ is finite, there is a choice of $t$ for which $\phi_t$ is injective on $\pi_j^{\oplus n}(B^{\oplus n})$ for every $j \in S \cup \{-1\}$.

Set $B' = \psi_t(B^{\oplus n})$, so that 
\[
\lvert \pi_j(B') \rvert = \lvert \phi_t \circ \pi_j^{\oplus n}(B^{\oplus n}) \rvert = \lvert \pi_j^{\oplus n}(B^{\oplus n}) \rvert
\]
and the result follows.
\end{proof}

\subsection{Equivalence of the \texorpdfstring{$n$}{n}-dimensional AKCs}
\begin{proof}[Proof of \autoref{thm:ndimAKC}]

We will use an argument from \cite{GreenRuzsa}. 
First, let $A\in \Z^n$ be a subset containing an $k$-term arithmetic progression with common difference $d$, for each $d \in \{1, \dots, N\}^n$.
Consider the map $f:\Z^n\to \Z$ defined by $f(x_1,\hdots,x_n)=\sum_{i=1}^n(10kN)^ix_i$.
This map creates a subset of $\Z$ that preserves the desired properties of $A$.
That is, $f(A)$ is a subset of $\Z$ containing an arithmetic progression for $N^n$ differences. 
Therefore, $|A|\geq F_k' (N^n)$, and thus $F_{n,k}(N)\geq F_k'(N^n).$ It follows that 
\[
\lim_{N\to \infty}\frac{\log F_{n,k}(N)}{\log N}\geq n\lim_{N\to \infty}\frac{\log F_k'(N^n)}{\log N^n},
\]
so that the $1$-dimensional AKC implies the $n$-dimensional AKC. 

For the converse, suppose we have a set $A_1\subset \Z$ containing a $k$-term progression with common difference $d$ for each $d\in \{1,\hdots,N\}$. Partition $\Z$ into intervals $I_j:=10kjN+\{1,\hdots,10kN\}$, for each $j\in \Z$. 
Then, each of these progressions is either entirely contained in some $I_j$, or split into two progressions where one is in $I_j$, the other in $I_{j+1}$, and one of these intersections has length $\geq k/2$. 
Let
\[
A_2=\bigcup_j(A_1 \cap I_j)-10kjN.
\]
$A_2$ is a subset of $\{1,\hdots,10kN\}$ containing a progression of length at least $k/2$ for each common difference $d\in\{1,\hdots,N\}$. Notice that $|A_2|\leq |A_1|.$

Let $M:=\lfloor{N^{1/n}}\rfloor$ and select $t\in \{-10kN,\hdots,20kN-1\}$ uniformly at random. Define 
\[
A_3(t):=\left\{(x_1,x_2, \hdots,x_n)\in\{0,\hdots M-1\}^n:\sum_{i=1}^n M^{i-1}x_i\in A_2+t\right\}.
\]
There are at least $(M/4k)^n$ such values of $d$ in $\{0,\hdots, N\}$ that can be written as $d=\sum_{i=1}^n M^{i-1}d_i$, with $0\leq d_i\leq M/2k$ for each $i$.
By assumption, a progression $\{x(d)+\lambda d:\lambda=0,1,\hdots \lfloor{k/2}\rfloor-1\}$ lies in $A_2$ for each such $d$. 
Thus, $\{x(d)+t+\lambda d:\lambda=0,1,\hdots \lfloor{k/2}\rfloor-1\}$ lies in $A_2+t$. Let 
\[
S:=\left\{\sum_{i=1}^n M^{i-1}s_i :\forall i,0\leq s_i <M/2  \right\}.
\]

If $t\in -x(d)+S$, then $A_3(t)$ contains a progression of length $k/2$ and difference $(d_1,\hdots,d_n)$. That is, $\{(s_1,\hdots s_n)+\lambda(d_1,\hdots,d_n):\lambda\in \{0,1,\hdots,k-1\}\}$, where $x(d)+t=\sum_{i=1}^n M^{i-1}s_i$.

Notice that $-x(d)+S\subset\{-10kN,\hdots,20kN-1\}$, because $0\leq x(d)\leq 10kN$ and $S\subset \{0,1,\hdots, M^n\}$. It follows that
\[
\Pr(t\in-x(d)+S)=\frac{1}{30kN}\lvert S\rvert\leq \frac{1}{30kN}(\frac{N}{2})^n \gg_{k,n} 1.
\]
Summing over the $(M/2k)^n\gg_{k,n} N$ choices of $d$, the expected number of $d$ for which $t\in-x(d)+S$ is $\gg_{k,n} N$. Fix some choice of t such that $t\in -x(d)+S$ for $\gg_{k,n} N\gg_{k,n} M^n$ values of $d$, and write $A_3:=A_3(t).$  Then, by construction, $|A_3|\leq|A_2|\leq|A_1|$, where $A_3$ contains a progression of length $\geq k/2$ and common difference $d$ for all $d$ in some set $D\subset \{0,\hdots M-1\}^n$ with cardinality $\gg_{k,n} M^n$. Hence,
$F'_{n,k/2}(N)\leq F_k(N)$, which concludes the proof.
\end{proof}
%\subsection{Random percolations}\todo{Needs to be improved}
%We recall the definition of a fractal percolation. Given a parameter $p\in(0,1)$ and a collection $A_n$ of (half-open) intervals in $\R$, we subdivide each interval into 2 halves, and let $A_{n+1}$ be a random collection of intervals in which each subdivided half-interval of $A_n$ is included with probability $p$. Beginning with the collection $A_0:=\{[0,1)\}$ of intervals, we recursively construct a sequence $A_0,A_1,\hdots$, and define
%\[A^{\text{perc}(p)}:= \overline{\bigcap_{n\in\mathbb{N}} \tilde{A_n}}\]
%where $\tilde{A}_n=\bigcup A_n$ is the union of the intervals in $A_n$. At each iteration $n$, we construct a random measure $\nu_n(B)=p^{-n}\lambda(B\cap \tilde{A}_n)$. The distribution of such measures is denoted $\nu_n^{\text{perc}(p)}$. We can now state Shmerkin's theorem:
%\begin{theorem}
%Suppose $\nu_n^{(i)}\sim \nu_n^{\text{perc}(p)}$ independently for $i\in[m]$ and let $\mu_n=\prod_{i=1}^m \nu_n^{(i)}$. Let $\Gamma$ be a collection of lines in $\R^m$ not contained in any coordinate hyperplanes $H^I=\{x\in\R^m:x_i=0\text{ for all }i\in I\}$ for any $I\subsetneq[m]$. If $s>1-1/m$, then
%\[ Y_n^V=\int_V \mu_n d\mathcal{H}^1 \]
%converges to $Y^V$ uniformly over all $V\in\Gamma$. Furthermore, there exists a nonempty open neighborhood $U\subset \R^m$ such that every $V\in\Gamma$ which meets $U$ satisfies $\mathbb{P}(Y^V>0)>0$.
%\end{theorem}

\bibliography{bib}{}
\bibliographystyle{plain}

\end{document}